%% file: LevelBetti.tex
\documentclass[12pt, a4paper]{amsart}

\usepackage{amscd}
\usepackage{a4}
\usepackage{amsmath}
\usepackage{amsxtra}
\usepackage[latin1]{inputenc}
\usepackage{graphicx}
\usepackage{amssymb}
\usepackage[all]{xy}
\usepackage{color}

\numberwithin{equation}{section}

\newtheorem{thm}{Theorem}[section]
\newtheorem{lemma}[thm]{Lemma}
\newtheorem{prop}[thm]{Proposition}
\newtheorem{cor}[thm]{Corollary}
\newtheorem{conj}[thm]{Conjecture}

\theoremstyle{definition}

\newtheorem{dfn}[thm]{Definition}

\newtheorem{rmk}[thm]{Remark}
\newtheorem{example}[thm]{Example}

\DeclareMathOperator{\Hom}{Hom}
\DeclareMathOperator{\HHom}{\sideset{^*}{}\Hom}

\DeclareMathOperator{\Tor}{Tor}

\DeclareMathOperator{\Ext}{Ext}

\newcommand{\MacRep}[4]{[{#1}_{(#2)}]_{#4}^{#3}}
\newcommand{\up}{\overline{d}}
\newcommand{\down}{\underline{d}}
\newcommand{\sspec}{_{\down, \up}}

\newcommand{\vecsp}{V\sspec}

\newcommand{\poset}{\varPi \sspec}
\newcommand{\pure}{\pi}

\newcommand{\pdim}{p}
\newcommand{\MaxLex}{\beta^{\mathrm{lex}}}

\begin{document}

\title{Graded Betti numbers and $h$-vectors of level modules} \date{\today}

\author{Jonas S\"{o}derberg} \address{Department of Mathematics, KTH
  \\ S--100 44 Stockholm \\ Sweden} \email{jonasso@math.kth.se}


\begin{abstract}
  We study $h$-vectors and graded Betti numbers of level modules up to
  multiplication by a rational number. Assuming a conjecture on the
  possible graded Betti numbers of Cohen-Macaulay modules we get a
  description of the possible $h$-vectors of level modules up to
  multiplication by a rational number. We also determine, again up to
  multiplication by a rational number, the cancellable $h$-vectors and
  the $h$-vectors of level modules with the weak Lefschetz property.
  Furthermore, we prove that level modules of codimension three
  satisfy the upper bound of the Multiplicity conjecture of Herzog,
  Huneke and Srinivasan, and that the lower bound holds if the module,
  in addition, has the weak Lefschetz property.
\end{abstract}

\maketitle


\section{Introduction} 
Level algebras, introduced by Stanley \cite{Stanley2}, often shows up in algebraic
geometry and combinatorics, and have in recent years received a lot of
attention themselves, especially the problem of describing their
$h$-vectors. Iarrobino \cite{Iarrobino-1} (see also Chipalkatti and
Geramita \cite{Chipalkatti-Geramita}) determined all $h$-vectors of
level algebras of codimension two, and this result was generalized to
level modules by the author \cite{Soderberg2}. In codimension three,
Stanley determined the $h$-vectors of Gorenstein algebras \cite{Stanley}, that is,
level algebras of type one. There are many other results in these
directions (see \cite{Geramita-Harima-Migliore-Shin} for an overview)
   but so far the problem of describing the $h$-vectors of all
level algebras seems out of reach.

In \cite{Boij-Soderberg}, Boij and the author gives a conjecture
(Conjecture \ref{MainConj} below) on the possible graded Betti numbers of
Cohen-Macaulay modules up to multiplication by a positive rational
number, and show that this conjecture implies the Multiplicity
conjecture of Herzog, Huneke and Srinivasan \cite{Herzog-Srinivasan},\cite{Huneke-Miller}. 
The heuristic is that the
set of possible graded Betti numbers, or $h$-vectors for that matter, is
easier to describe, when we consider not only algebras but also
modules, and care only for a description up to multiplication by a
positive rational number. For example, the $h$-vectors of
Cohen-Macaulay modules, generated in degree zero, of
codimension $p$, are, up to multiplication by a positive rational
number, the finite sequences $\{h_i\}_{i \in \mathbb{Z}}$ satisfying
$$ 
 \frac{h_i}{r_{i}} - \frac{h_{i+1}}{r_{i+1}} \geq 0
$$ 
for all $i$, where $r_i$ is the dimension, as a vector space, of the
graded component of degree $i$ of the polynomial ring in $p$
variables.  This should be compared with Macaulay's description \cite{Macaulay} (and
Hulett's generalization to modules of this description and \cite{Hulett-2}) which
precisely characterizes these sequences but which is more
combinatorial in nature and more complicated.

In this paper we take as our starting point the conjecture on the
possible graded Betti numbers of
Boij and the author, and ask what can be said about the
$h$-vectors and graded Betti numbers of level modules if this conjectured
is assumed to be true. This does not mean that all our results depend
on this conjecture, but rather that we are guided by it. In particular,
we look for descriptions of $h$-vectors and graded Betti numbers
only up to multiplication by a rational number.

In Section 4 we determine the $h$-vectors of level modules
which satisfy the condition of the conjecture
of Boij and the author. It turns out that, up to multiplication with a positive
rational number, these $h$-vectors are precisely those which are
non-negative linear combinations of $h$-vectors of extremely compressed
level modules of the same socle degree, and moreover, this condition
can be expressed as a set of linear inequalities, each in three consecutive
entries of the $h$-vector (Theorem \ref{NormHVec}).  We also show
that the Betti numbers of the linear combinations of extremely
compressed level modules, bounds from above the Betti numbers of the
level module in this case and even that the Betti numbers of the level
module are obtained from these by a sequence of consecutive
cancellations (Proposition~\ref{MaxCompressedProp}). 
Unfortunately we can not prove that level modules in general satisfy this
condition, but the connection with extremely compressed level modules
is interesting, and it might be easier prove this result for level
modules in general than to prove the conjecture of Boij and the author.

When considering graded Betti numbers up to multiplication
by a rational number, the maximal
ones, given the $h$-vector, have a simple description. This observation 
turns out to be very useful and is
used in almost all parts of this article.
 
With this description of maximal Betti numbers in mind, the notion of
\emph{cancellable} $h$-vectors, introduced by Geramita and Lorenzini
\cite{Geramita-Lorenzini} and also studied by the author in
\cite{Soderberg1}, is revisited and we give a description, this time
up to multiplication by a rational number, of all cancellable
$h$-vectors (Theorem~\ref{DualCancelProp}).

In section 6 we study the $h$-vectors, and graded Betti numbers,
of artinian level modules with 
the weak Lefschetz property.
In \cite{Harima-Migliore-Nagel-Watanabe} the possible $h$-vectors
of algebras with the weak and strong Lefschetz property are
determined and a sharp upper bound on the graded Betti numbers
of these algebras is given. The possible $h$-vectors are also
known for Gorenstein algebras with the weak Lefschetz property \cite{Harima}, 
but for level algebras in
general there are very few results. 
We determine the possible $h$-vectors of artinian level modules with
the weak Lefschetz property up to multiplication by a rational number,
and furthermore, we give an upper bound on the graded Betti numbers of these modules 
(Theorem \ref{WLPUpToMultProp}).
This upper bound is sharp in the sense that it is a rational multiple of the Betti diagram of
some level module with the weak Lefschetz property, and furthermore, the upper bound
satisfy the condition of the conjecture of Boij and the author.
It turns out that the existence of this upper bound is enough to prove
the conjecture of Boij and the author for level modules of codimension three with
the weak Lefschetz property which in turn proves the Multiplicity conjecture
for these modules (Proposition \ref{CodimThreeWLPCor}).

Finally, in Section 7 we restrict our attention to level modules of codimension
three.  We observe that the conjecture of Boij and the author
implies a strengthening of the Multiplicity
conjecture proposed by Zanello in the case of level algebras of
codimension three \cite{Migliore-Nagel-Zanello}.
We finish by proving that the upper bound of the
Multiplicity conjecture holds for any level module of codimension three
(Theorem~\ref{MCUpperThm}).

\section{Preliminaries} \label{Preliminaries} We begin with the basic
definitions.  
Let $R = k[x_1, x_2, \dots, x_n]$ be the polynomial ring
in $n$ variables over a field $k$. Consider $R$ as a graded ring by
giving each $x_i$ degree one and let 
$\mathfrak{m} = (x_1,x_2,\dots,x_n)$ be the unique graded maximal
ideal. All $R$-modules in this article are
assumed to be finitely generated and graded.  The $d$-th twist of an
$R$-module $M$, denoted by $M(d)$, is defined by $M(d)_i = M_{i+d}$.
If the $R$-module $M$ has a minimal free resolution given by
$$
0 \to \bigoplus_j R(-j)^{\beta_{\pdim,j}} \to \dots \to \bigoplus_j
R(-j)^{\beta_{0,j}} \to M \to 0
$$
then $\beta^R_{i,j}(M) = \beta_{i,j}$ are the graded Betti numbers of
$M$. When the ring in consideration is clear from the context we
simply omit the superscript and write $\beta_{i,j}(M)$ for the Betti
numbers of $M$.  From the minimal free resolution of $M$ we also
see that the projective dimension of $M$ is $\pdim$ and when $M$ is
Cohen-Macaulay it is equal to the codimension, that is, $\dim R - \dim
M = \pdim$. When considering artinian $R$-modules of codimension
$p$ we consequently have $n=p$. Furthermore, if $M$ is Cohen-Macaulay then there is an element
$h_M(t) \in \mathbb{Z}[t,t^{-1}]$ such that
$$
h_M(t) = \frac{S_M(t)}{(1-t)^\pdim},
$$
where $S_M(t) = \sum_{i,j} (-1)^i \beta_{i,j}(M) t^j$.  The element
$h_M(t)$ is called the $h$-vector of $M$ and
$$
e(M) = h_M(1)
$$
is called the multiplicity of $M$. The matrix $\beta(M)$ with entries
$\beta(M)_{i,j} = \beta_{i,j}(M)$ for each integer $i$ and $j$, is
called the Betti diagram of $M$. 

The Hilbert function of $M$ is the function $H(M,\_): \mathbb{Z} \to \mathbb{Z}$
defined by $H(M,d) = \dim_k M_d$. If $M$ is an artinian $R$-module
with $h$-vector $h = \sum_{i\in \mathbb{Z}}h_it^i$, then $H(M,d) = h_d$.

Let $M$ be a Cohen-Macaulay $R$-module of codimension $p$. We define 
the dual of $M$, denoted by $M^{\vee}$, to be the $R$-module $\Ext_R^p(M,R)$.
The dual of of $M$ is in fact Cohen-Macaulay and its Betti diagram is
obtained from that of $M$ by $\beta_{i,j}(M^{\vee}) =
\beta_{p-i,p-j}(M)$.

The \emph{maximal and minimal shifts} of degree $i$ of $M$ are defined
by $\up_i(M) = \max\{j\,|\, \beta_{i,j}(M) \neq 0\}$ and $\down_i(M) =
\min\{j\,|\, \beta_{i,j}(M) \neq 0\}$, respectively. It is well known
that when $M$ is Cohen-Macaulay then both of the sequences $\up(M) =
(\up_0(M),\up_1(M),\dots,\up_\pdim(M))$ and $\down(M) =
(\down_0(M),\down_1(M),\dots,\down_\pdim(M))$ are strictly increasing.

\section{Pure diagrams and consecutive cancellations}
\label{PureDiagrams} We will now explain a conjecture of Boij and the
author on the set of possible Betti diagrams of Cohen-Macaulay modules
up to multiplication by a rational number.  We also explain how Betti
diagrams are obtained from, entry by entry, larger ones with the same
$h$-vector by consecutive cancellations.  We refer to
\cite{Boij-Soderberg} for a more detailed exposition of these matters.

We begin by explaining the notion of a diagram.  Let $D$ be a matrix
with rational entries satisfying $D_{i,j} = 0$ when $i<0$ or $i>p$ and
assume furthermore that $D$ has only a finite number of non-zero
entries.  Denote by $S_D(t)$ the element $S_D(t) = \sum_{i,j} (-1)^i
D_{i,j} t^j$ in $\mathbb{Q}[t,t^{-1}]$.  If $S_D(t)$ is divisible by
$(1-t)^p$ we call $D$ a \emph{diagram} of codimension $p$. This
definition is motivated by the fact that the Betti diagram of a
Cohen-Macaulay module of codimension $p$ is a diagram of codimension
$p$. Let $\up = (\up_0,\up_1,\dots,\up_p)$ and $\down =
(\down_0,\down_1,\dots,\down_p)$ be two strictly increasing sequences
of integers and denote by $\vecsp$ the vector space over the rational
numbers of all diagrams $D$, of codimension $p$, satisfying $D_{i,j} =
0$ whenever $j < \down_i$ or $j>\up_i$.

An $R$-module $M$ has a pure resolution of type $d =
(d_0,d_1,\dots,d_p)$ if its minimal free resolution has the form
$$
0 \to R(-d_\pdim)^{\beta_{\pdim}} \to \dots \to R(-d_0)^{\beta_{0}}
\to M \to 0,
$$
that is, if the maximal and minimal shifts of $M$ are equal.

If $D$ is a diagram whose maximal and minimal shifts both are equal to
the strictly increasing sequence of integers $d =
(d_0,d_1,\dots,d_p)$, then we call $D$ a pure diagram of type $d$.
The dimension of the vector space $\vecsp$, for any strictly
increasing sequences $\up$ and $\down$, is $1+\sum_{i=0}^\pdim
(\up_i-\down_i)$ (see \cite[Prop. 2.7]{Boij-Soderberg}), so the
dimension of the vector space of all pure diagrams of type $d$,
$V_{d,d}$, is one. Hence any two pure diagrams of the same type are
rational multiples of each other, and we denote by $\pure(d)$ the pure
diagram of type $d$ satisfying $\pure(d)_{0,d_0} = 1$.

Boij and the author have the following conjecture in
\cite{Boij-Soderberg}.

\begin{conj}[Boij, S\"{o}derberg]\label{MainConj}
  The Betti diagram of any Cohen-Macaulay $R$-module is a non-negative
  linear combination of pure diagrams and furthermore, any pure
  diagram is a rational multiple of the Betti diagram of some
  Cohen-Macaulay $R$-module.
\end{conj}
 
Let $M = F/N$ be a Cohen-Macaulay $R$-module, where $F$ is a free
$R$-module and $N$ a submodule of $F$.  Then there is a lexicographic
submodule $L$ of $F$ such that $F/L$ and $M$ have the same $h$-vector
(see \cite{Macaulay} and \cite{Hulett-2}).  The Betti diagram of $F/L$
is completely determined by the $h$-vector, $h(t)$, of $M$ and of the
free module $F$ and we denote this diagram by $\MaxLex_F(h(t))$.  The
Betti numbers of $M$ are smaller than the Betti numbers of $F/L$, and
we even have that $\beta(M)$ is obtained from $\MaxLex_F(h(t))$ by a
sequence of consecutive cancellations (see \cite{Peeva}).  A
consecutive cancellation is defined as follows.  Let $k$ and $l$ be
two integers and define a diagram $C^{k,l}$ by
$$
(C^{k,l})_{i,j} = \begin{cases} 1 &\text{ when $(i,j) = (k,l)$ or
    $(i,j)=(k+1,l)$}, \\ 0 &\text{ otherwise}.\end{cases}
$$
A consecutive cancellation in position $(k,l)$ of a diagram $D$ is a
diagram
$$
D^{\prime} = D - bC^{k,l},
$$
where $b$ is a non-negative rational number, and we assume that
$D^\prime$ have no negative entries.  Note that a consecutive
cancellation does not change the polynomial $S_{D}(t) = \sum_{i,j}
(-1)^i D_{i,j} t^j$, that is, $S_{D}(t) = S_{D^{\prime}}(t)$, and
hence not the $h$-vector or multiplicity of $D$. In short, for any
Cohen-Macaulay $R$-module $M$ with $h$-vector $h(t)$ we have that
\begin{equation} \label{MaxLexBettiEq} \beta(M) = \MaxLex_F(h(t)) -
  \sum_{i,j} b_{i.j} C^{i,j}
\end{equation}
for some non-negative integers $b_{i,j}$.

Let $F$ be a free $R$-module with basis $e_1,e_2,\dots,e_t$ and assume
that the basis of $F$ is ordered by $e_1 < e_2 < \dots < e_t$ where
$\deg e_1 \leq \deg e_2 \leq \dots \leq \deg e_t$.  Let
$H:\,\mathbb{Z} \to \mathbb{Z}$ be the Hilbert function of some
submodule of $F$. Then, as mentioned above, there is a lexicographic
submodule $L$ with $H$ as Hilbert function and we will now describe
$L$. The monomials of $F$ are the elements on the form $ue_i$ for some
monomial $u$ of $R$ and basis element $e_i$ of $F$.  The lexicographic
order on monomials of $F$ is given by $ue_i > ve_j$ when $i<j$ or when
$i=j$ and $u>v$ in the lexicographic order on monomials of $R$.  The
lexicographic submodule $L$ is generated in each degree $d$ by the
$H(d)$ largest monomials in $F_d$.

For any degree $d$, there are unique integers $k$ and $q$ such that $0
\leq q < H(R,d-\deg e_k)$ and
\begin{equation} \label{LexEq} H(d) = \sum_{i=1}^{k-1} H(R,d-\deg e_i)
  + q.
\end{equation}
For any $F$, $H$ and $d$ denote these numbers by $k(F,H,d)=k$ and
$q(F,H,d) = q$ and furthermore let $g(F,H,d) = \deg e_k$.  This means
that $L$ is generated in degree $d$ by
$$
\sum_{i=1}^{k-1} \mathfrak{m}^{d-\deg e_i}e_i + Ie_{k}
$$
where $k = k(F,H,d)$ and $I$ is the ideal generated by the $q(F,H,d)$
largest monomials of degree $d-g(F,H,d)$ in $R$.

In what follows, let $e_{i,j}$, where $i=1,2,\dots,t$, $j =
1,2,\dots,s$ and $\deg e_{i,j} = \deg e_i$, be a basis of $F^s$
ordered by $e_{i,j} > e_{i',j'}$ when $i < i'$ or when $i = i'$ and $j
> j'$, for any integer $s$.

\begin{lemma} \label{LexLemma} Let $N$ be a submodule of $F$ such that
  the module $M = F/N$ is artinian. Then there is an integer $s$ and a
  lexicographic submodule $L$ of $F^s$ such that $L$ and $N^s$ have
  the same $h$-vector, and $L$ is on the form
$$
L = \sum_{j=1}^s\sum_{i=1}^t \mathfrak{m}^{c_{i,j}}e_{i.j}
$$
where $c_{i,j}$ are non-negative integers.
\end{lemma}

\begin{proof}
  Let $H$ be the Hilbert function of $M$.  To prove the lemma it is
  enough to show that there is a positive integer $s$ such that
  $q(F^s,sH,d) = 0$ for all non-negative integers $d$.  Note that in
  (\ref{LexEq}), $H(d) = H(F,d)$ implies $q=0$, and hence that $H(d)=
  H(F,d)$ implies $q(F,H,d) = 0$. Since $M$ is artinian, $H(d) \neq
  H(F,d)$ only for a finite number of integers and hence the same
  holds for $q(F,H,d)$. This means that there exists a positive
  integer $s$ such that
  \begin{equation} \label{sEq} s\frac{q(F,H,d)}{H(S,d-g(F,H,d))}
  \end{equation}
  is an integer for all $d$. Let $s$ be such an integer.

  Fix a degree $d$ and let $k = k(F,H,d)$ and $q = q(F,H,d)$ then we
  get by multiplying both sides of (\ref{LexEq}) by $s$,
$$
sH(d) = s\sum_{i=1}^{k-1} H(R,d-\deg e_i) + sq =
\sum_{j=1}^s\sum_{i=1}^{k-1} H(R,d-\deg e_{i,j}) + sq
$$
and since by (\ref{sEq}), $sq = s'H(R,d-\deg e_k)$ for some $0 \leq s'
< s$ we get
$$
sH(d) = \sum_{j=1}^s\sum_{i=1}^{k-1} H(R,d-\deg e_{i,j}) +
\sum_{i=j}^{s'}H(R,d-\deg e_{k,j}).
$$
From the above equation it follows that $q(F^s,sH,d) = 0$ for all
non-negative integers $d$. This means that $L$ is generated in degree
$d$ by
$$ \sum_{j=1}^s\sum_{i=1}^{k-1} \mathfrak{m}^{d-\deg e_{i,j}}e_{i,j} +
\sum_{i=j}^{s'} \mathfrak{m}^{d-\deg e_{k,j}}e_{i,j}.
$$
and since this holds for any degree $d$ the lemma follows.
\end{proof}

\begin{prop}\label{NormMaxBettiProp}
  Let $M = F/N$ be an artinian $R$-module of codimension $p$ with
  $h$-vector $h(t)$.  Then there exists non-negative rational numbers
  $a_{i,j}$ and an integer $s$ such that the diagram
$$
D = \sum_{i,j} a_{i,j} \beta(R/\mathfrak{m}^{j+1}(-\deg e_i))
$$
have $h$-vector $h(t)$ and
\begin{equation*}
  \frac{1}{s}\MaxLex_{F^s}(sh(t)) = D.
\end{equation*}
As a consequence we get that $\beta(M)$ is obtained from $D$ by a
sequence of consecutive cancellations, that is,
$$
\beta(M) = D - \sum_{i,j}b_{i,j}C^{i,j}
$$
for some non-negative rational numbers $b_{i,j}$.
\end{prop}

\begin{proof}
  By Lemma~\ref{LexLemma} there is an integer $s$ and a lexicographic
  submodule $L$ of $F^s$ such that $F^s/L$ have $h$-vector $sh(t)$ and
  $L$ is on the form
$$
L = \sum_{j=1}^s\sum_{i=1}^t \mathfrak{m}^{c_{i,j}}e_{i,j}
$$ 
where $c_{i,j}$ are non-negative integers.  We get
$$
\MaxLex_{F^s}(sh(t)) = \beta(F^s/L)
$$
and since then
$$
F^s/L =
\bigoplus_{j=1}^s\bigoplus_{i=1}^tR/\mathfrak{m}^{c_{i,j}}(-\deg
e_{i,j})
$$
we see that
$$
\MaxLex_{F^s}(sh(t)) = \sum_{j=1}^s\sum_{i=1}^t
R/\mathfrak{m}^{c_{i,j}}(-\deg e_{i,j})
$$
Collecting all the terms in the sum above where the $c_{i,j}$ and
$\deg e_{i,j}$ are the same, yields
$$
\MaxLex_{F^s}(sh(t)) = \sum_{i,j}a'_{i,j} R/\mathfrak{m}^{j+1}(-\deg
e_{i})
$$
for some non-negative integers $a'_{i,j}$.  Now, let $a_{i,j} =
a'_{i,j}/s$ and $D = \MaxLex_{F^s}(sh(t))/s$.

Since the Hilbert function of $M^s$ is $sh(t)$ we get that
$\beta(M^s)$ is obtained from $\MaxLex_{F^s}(sh(t)) = sD$ by a
sequence of consecutive cancellations, that is,
$$
\beta(M^s) = sD - \sum_{i,j}b'_{i,j}C^{i,j}
$$
for some non-negative integers $b_{i,j}$. Now $\beta(M^s) = s\beta(M)$
so dividing both sides of the equation above with $s$ shows that
$\beta(M)$ is obtained from $D$ by a sequence of consecutive
cancellations.
\end{proof}

When $F$ is generated in degree zero we omit the subscript $F$ and
simply write
$$
\MaxLex(h(t)) = \MaxLex_F(h(t))
$$
and moreover, in this case there is a simple expression for the
numbers $a_{i,j}$ of Proposition~\ref{NormMaxBettiProp}.

\begin{prop}\label{SingleMaxBettiProp}
  Assume that $F$ is generated in degree zero and let $M=F/M$ be an
  artinian $R$-module, of codimension $p$, with $h$-vector the
  polynomial $h(t) = \sum_{i \in \mathbb{Z}} h_it^i$ of degree $c$.
  Then
  \begin{equation*}
    \frac{1}{s}\MaxLex(sh(t)) = \sum_{j=0}^c a_j \beta(R/\mathfrak{m}^{j+1})
  \end{equation*}
  for some integer $s$, if and only if
$$
a_j = \frac{h_j}{r_j} - \frac{h_{j+1}}{r_{j+1}}
$$
where $r_j = \binom{p-1+j}{p-1}$. As a consequence we get that the
polynomial $h(t)$ is a rational multiple of an $h$-vector of an
artinian $R$-module, generated in degree zero, if and only if the
$a_j$'s above are non-negative.
\end{prop}

\begin{proof}
  By Proposition~\ref{NormMaxBettiProp} there is an integer $s$ and a
  diagram
$$
D = \sum_{j} a_j \beta(R/\mathfrak{m}^{j+1}),
$$
for some non-negative rational numbers $a_j$, such that
$$
\frac{1}{s}\MaxLex(sh(t)) = D
$$
Denote by $h_j$ the $h$-vector of $R/\mathfrak{m}^{j+1}$ and observe
that it is given by
$$
h_{j}(t) = r_0 + r_1t + \dots + r_jt^j.
$$    
The $h$-vector of $D$ is then
$$
h_D(t) = \sum_{j}a_jh_{j}(t).
$$
Solving the linear equation $h(t) = \sum_{j}a_jh_{\pi_j}(t)$ for the
rational numbers $a_j$ gives the unique solution $a_j = h_j/r_j -
h_{j+1}/r_{j+1}$.

We also see that the $a_j$'s are non-negative if $h(t)$ is the
$h$-vector of an artinian $R$-module generated in degree zero, and
furthermore, if the $a_j$'s are non-negative then the module
$$
N = \bigoplus_{j=0}^c \bigoplus_{i=1}^{ma_j} R/\mathfrak{m}^{j+1}
$$
where $m$ is an integer such that $ma_j$ is an integer for
$j=0,1,\dots,c$, have $h$-vector $h_N(t) = mh(t)$ which shows that $h$
is a rational multiple of the $h$-vector of an artinian $R$-module.
\end{proof}


\section{$h$-vectors of level modules}
An $R$-module $M$ is \emph{level} if it is Cohen-Macaulay and its
artinian reduction has socle and generators concentrated in single
degrees. This can be expressed in terms of the Betti diagram of $M$.
If $\down = (\down_0,\down_1,\dots,\down_p)$ and $\up =
(\up_0,\up_1,\dots,\up_p)$ are the minimal and maximal shifts of $M$,
then $M$ is level if and only if $\down_0 = \up_0$ and $\down_\pdim =
\up_\pdim$.

We will now show that if $M$ is level and its Betti diagram is a
non-negative linear combination of pure diagrams, which, according to
Conjecture~\ref{MainConj}, every Betti diagram of a Cohen-Macaulay
module is, then its $h$-vector satisfies a certain condition and any
polynomial satisfying this condition is, after multiplication by some
integer, the $h$-vector of a level $R$-module.  This condition on the
$h$-vector has to do with extremely compressed level modules, so we
begin by explaining this notion.

A compressed level module is an artinian level module of maximal
Hilbert function
among all artinian level modules with a given codimension, number of
generators, socle type and socle degree.  (The socle type is the
dimension of the socle as a $k$-vector space and the socle degree, the
degree in which the socle is concentrated).  It turns out that this
Hilbert function is given by $i \mapsto \min \{sr_i, tr_{c-i}\}$ where
$s$ is the number of generators, $t$ is the socle type and $c$ is the socle
degree.  For any $s$, $t$ and $c$ there is a level module with this
Hilbert function \cite{Boij}, and, as already mentioned, this module
is said to be compressed.  If in addition $sr_i = tr_{c-i}$ for some
$i$ then the level module having this Hilbert function is said to be
\emph{extremely compressed}.

\begin{prop} \label{CompExistProp} For any integers $d_0$, $j$ and
  $c$, such that $d_0 \leq j \leq c$, there is a Cohen-Macaulay
  $R$-module with pure resolution of type $d =
  (d_0,j+1,j+2,\dots,j+p-1,c+p)$ and the artinian reduction of any
  Cohen-Macaulay $R$-module with a pure resolution of this type is an
  extremely compressed level module.
\end{prop}

This follows from the work on compressed level modules in \cite{Boij}.
We include a proof here for the reader's convenience.

\begin{proof}
  By artinian reduction we may assume that $p=n$, and hence only
  consider artinian $R$-modules and by shifting the degrees we may
  assume that $d_0 = 0$. Let $M$ be the compressed level $R$-module
  with socle degree $c$, $s$ number of generators and socle type $t$,
  where $s = r_{c-j}$ and $t = r_j$.  Then $sr_j = tr_{c-j}$ so $M$ is
  extremely compressed.  To see that the minimal free resolution of
  $M$ is pure of type $d= (0,j+1,j+2,\dots,j+p-1,c+p)$, consider
  $\down_1(N)$ and $\up_{p-1}(N)$. Note that $\down_1(N)$ is the
  lowest degree of a relation among the generators of $M$ and that the
  Hilbert function of $M$ is that of a free module in degrees lower
  than, or equal to, $j$, where $j$ is the integer such that $sr_j =
  tr_{c-j}$.  This implies that $\down_1(N) = j+1$. Since the Betti
  numbers in column $p-1$, of the Betti diagram of $M$, describes the
  degrees of the relations among the generators of the dual of $M$, we
  can use the same argument again and we get $\up_{p-1}(N) = p+j-1$.
  Now there is only one possibility for the minimal and maximal shifts
  of $M$, since they are both strictly increasing sequences of
  integer, they must be the same and equal $d =
  (0,j+1,j+2,\dots,j+p-1,c+p)$.
\end{proof}

\begin{dfn}
For any integers $d_0 \leq j \leq c$ we call a pure diagram of type 
$d = (d_0,j+1,j+2,\dots,j+p-1,c+p)$ an \emph{extremely compressed} diagram.
\end{dfn}

We will now prove a proposition from which the promised description of the
$h$-vectors of level modules, up to multiplication by a positive
rational number, follows easily, as we will soon see.

\begin{prop}\label{MaxCompressedProp}
  Let $D$ be a non-negative linear combination of pure diagrams of
  codimension $p$. Then there is a diagram $E$ which is a non-negative
  linear combination of extremely compressed diagrams of
  codimension $p$ such that $D$ is obtained from $E$ by a sequence of
  consecutive cancellations and furthermore, $E_{0,j} = D_{0,j}$ and $E_{p,j} =
  D_{p,j}$ for all integers $j$.
\end{prop}

We will need the following observation.
\begin{lemma} \label{CompLemma} Let $\pure(d)$ be an extremely
  compressed diagram of type
  $d = (d_0,j+1,j+2,\dots,j+p-1,d_p)$.  Then we can write $\pure(d)$ on the form
$$
\pure(d) = \sum_{j} a_j \pure(d_0,j+1,j+2,\dots,j+p-1,j+p) -
 \sum_{i,j}b_{i,j}C^{i,j}
$$
where $a_i$ and $b_{i,j}$ are non-negative
rational numbers.
\end{lemma}

\begin{proof}
  By Proposition~\ref{CompExistProp} there exists an extremely
  compressed level module, $M$, with a pure resolution of type $d$, and hence
  $\beta(M) = q\pure(d)$ for some positive integer $q$.  By
  Proposition~\ref{NormMaxBettiProp} we have
$$
\beta(M) = \sum_{j} a'_{j}\beta(R/\mathfrak{m}^{j+1}(-d_0))
- \sum_{i,j}b_{i,j}C^{i,j}
$$
for some non-negative rational numbers $a'_i$ and $b_{i,j}$.
 The lemma follows since $\pure(d) = \beta(M)/q$ and the
Betti diagram $\beta(R/\mathfrak{m}^{j+1}(-d_0))$ is
the pure diagram
$$
\pure(d_0,d_0 + j+1,d_0+j+2,\dots,d_0+j+p-1,d_0+j+p)
$$
for each integer $j$.
\end{proof}

For the induction step in the proof of
Proposition~\ref{MaxCompressedProp}, and for several other things
later, we need a result from \cite{Boij-Soderberg}.

\begin{dfn} \label{stepDown} For any sequence of integers
  $d=(d_0,d_1,\dots,d_{k-1},d_k,d_{k+1},\dots,d_p)$ denote by
  $\tau_k(d)$ the sequence
  $(d_0,d_1,\dots,d_{k-1},d_{k+1},\dots,d_p)$.  If $\up$ and $\down$
  are such that $\up_k = \down_k$ for some integer $k$, then there is
  an isomorphism of vector spaces (see \cite[Lemma
  3.3]{Boij-Soderberg})
$$
\phi_k:\, \vecsp \to V_{\tau_k(\down),\tau_k(\up)}
$$
given by
$$ 
\phi_k(D)_{i,j} = |\up_k-j| \cdot
\begin{cases} D_{i,j} &\text{ when $i<k$} \\
  D_{i+1,j}& \text{ otherwise} \end{cases}
$$
for any $D \in \vecsp$.  Note that the image of a pure diagram under
this map is $\phi_k(\pure(d)) = \up_k\pure(\tau_k(d))$, which is a
pure diagram of codimension one lower than the pure diagram we started
with.
\end{dfn}

\begin{proof}[Proof of Proposition~\ref{MaxCompressedProp}]
  It is enough to prove the case where $D = \pure(d)$ for some pure
  diagram of type $d=(d_0,d_1,\dots,d_p)$.  We use induction on the
  codimension $p$. If $p = 0,1$ or $2$ then any pure diagram is
  extremely compressed and the assertion follows. Assume now that $p >
  2$. Let $\phi_p$ be the linear map from Definition~\ref{stepDown}.
  By induction we may assume that the assertion of the proposition
  holds for the codimension $p-1$ pure diagram
  $\pure(d_0,d_1,\dots,d_{p-1})$ and hence that $\phi_p(\pure(d)) =
  d_p\pure(d_0,d_1,\dots,d_{p-1})$ is obtained from a non-negative
  linear combination of extremely compressed diagrams by a
  sequence of consecutive cancellations, that is,
  \begin{equation} \label{sumOfCompEq} \phi_p(\pure(d)) =
    \sum_{j} a_j
    \pure(d_0,j+1,j+2,\dots,j+p-2,d_{p-1}) - \sum_{i,j}b_{i,j}C^{i,j}
  \end{equation}
  for some non-negative integers $a_i$ and $b_{i,j}$. Note that all
  these extremely pure diagrams have $d_0$ as their first shift and
  $d_{p-1}$ as their last, since otherwise we could never have that
  their linear combination would equal $d_p\pure(d_0,d_1,\dots,d_{p-1})$
  in position $(0,j)$ and $(p-1,j)$ for all $j$ and this holds by the
  induction assumption.  By applying Lemma~\ref{CompLemma} to each of
  the extremely compressed diagrams in (\ref{sumOfCompEq}), and
collecting terms with pure diagrams of the same type, we get
\begin{equation} \label{PhiPureEq}
\phi_p(\pure(d)) = \sum_{j} a'_j
 \pure(d_0,j+1,j+2,\dots,j+p-1) - \sum_{i,j}b'_{i,j}C^{i,j}
\end{equation}
where and $a'_i$ and $b'_{i,j}$ are non-negative
rational numbers.  Consider the inverse of $\phi_p$ and note that
$$
\phi_p^{-1}(\pure(d_0,j+1,j+2,\dots,j+p-1)) =
 \frac{1}{d_p}\pure(d_0,j+1,j+2,\dots,j+p-1,d_p)
$$
Note also that
$$
\phi_p^{-1}(C^{i,j}) = \frac{1}{d_p-j}C^{i,j}.
$$
Hence, applying $\phi_p^{-1}$ to both sides of (\ref{PhiPureEq}) gives
$$
\pure(d) = \phi_p^{-1}\phi_p(\pure(d)) = E -
\sum_{i,j}b''_{i,j}C^{i,j}
$$
where
$$
E = \sum_{j} a''_j
\pure(d_0,j+1,j+2,\dots,j+p-1,d_p)
$$ 
and $a''_j = a'_j/d_p$ and $b''_{i,j} = b'_{i,j}/(d_p-j)$.
We have that $E$ is a non-negative linear combination of extremely
compressed diagrams and that $\pure(d)$ is obtained from $E$ by a 
sequence of consecutive cancellations. It remains to prove
that $E_{0,j} = \pure(d)_{0,j}$ and $E_{p,j} = \pure(d)_{p,j}$
for all integers $j$. Note that $E_{0,j} = 0$ for
all $j\neq d_0$ and $E_{p,j} = 0$ for all $j \neq d_p$.
So the only entries of $E$ to consider are the one in position
$(0,d_0)$ and $(p,d_p)$.
To cancel the entry in position
$(0,d_0)$ with a consecutive cancellation we need a positive
entry in position $(0,d_0-1)$, and by considering the
shifts of the pure diagrams in the sum giving $E$ we see
that $E_{0,d_0-1} = 0$. Hence $E_{0,d_0} = \pure(d)_{0,d_0}$
which shows that $E_{0,j} = \pure(d)_{0,j}$
for all integers $j$. That  $E_{p,j} = \pure(d)_{p,j}$
for all integers $j$ follows in the same way.
\end{proof}

Since a cancellation of a Betti diagram does not affect its $h$-vector
we see that Proposition~\ref{MaxCompressedProp} implies that the
$h$-vector of a module, $M$, whose Betti diagram is a non-negative
linear combination of pure diagrams is a non-negative linear
combination of $h$-vectors of extremely compressed level modules. If
in addition, $M$ is level, all of these extremely compressed level
modules must be generated in the same degree and have the same socle degree. 
We will now describe the
$h$-vector of such a level module, and we assume for simplicity that
it is generated in degree zero.

The $h$-vector of the extremely compressed diagram $\pi_j =
\pure(0,j+1,j+2,\dots,j+p-1,c+p)$ is
\[
h_{\pi_j}(t) = r_0 + r_1t + \dots + r_jt^j +
\frac{r_j}{r_{c-j}}\left(r_{c-j-1}t^{j+1} + r_{c-j-2}t^{j+2} + \dots +
  r_0t^c \right).
\]
In fact, let $M$ be the extremely compressed level $R$-module with a
pure resolution of type $d = (0,j+1,j+2,\dots,j+p-1,c+p)$ from the
proof of Proposition \ref{CompExistProp}.  The $h$-vector of this
module is $h_N(t) = \sum_{i=0}^c h_it^i$ where $h_i =
\min\{r_{c-j}r_i,r_jr_{c-i}\}$ and since $h_{\pi_j}(t) =
h_N(t)/r_{c-j}$ we get the desired expression for $h_{\pi_j}(t)$.

\begin{lemma} \label{ECompLemma} For any polynomial $h(t) = h_0 + h_1t
  + \dots + h_ct^c$ we have that $h(t) = \sum_{i=0}^c f_ih_{\pi_i}(t)$
  for some rational numbers $f_0,f_1,\dots,f_c$, where $h_{\pi_i}(t)$
  is the $h$-vector of the extremely compressed diagram $\pi_i =
  (0,i+1,i+2,\dots,i+p-1,c+p)$, if and only if
$$
f_i = r_{c-i}\frac{\left|\begin{matrix}
      h_{i-1} & h_{i} & h_{i+1} \\
      r_{i-1} & r_i & r_{i+1} \\
      r_{c-i+1} & r_{c-i} & r_{c-i-1}
    \end{matrix}\right|}
{\left|\begin{matrix}
      r_{i-1} & r_i \\
      r_{c-i+1} & r_{c-i}
    \end{matrix}\right|
  \left|\begin{matrix}
      r_{i} & r_{i+1} \\
      r_{c-i-1} & r_{c-i}
    \end{matrix}\right|}.
$$
\end{lemma}

\begin{proof}
  This is just a question of solving a system of linear equations.
  Let $h_{\pi_j,i}$ be the coefficient of $t^i$ in the polynomial
  $h_{\pi_j}(t)$, for each $j$, and consider the determinant
$$
\left|\begin{matrix}
    h_{i-1} & h_{i} & h_{i+1} \\
    r_{i-1} & r_i & r_{i+1} \\
    r_{c-i+1} & r_{c-i} & r_{c-i-1}
  \end{matrix}\right|
= \sum_{j=0}^c f_j \left|\begin{matrix}
    h_{\pi_j,i-1} & h_{\pi_j,i} & h_{\pi_j,i+1} \\
    r_{i-1} & r_i & r_{i+1} \\
    r_{c-i+1} & r_{c-i} & r_{c-i-1}
  \end{matrix}\right|.
$$
Whenever $i \neq j$ the vector $(h_{\pi_j,i-1},h_{\pi_j,i},
h_{\pi_j,i+1})$ is a multiple of either $(r_{i-1},r_i,r_{i+1})$ or
$(r_{c-i+1},r_{c-i},r_{c-i-1})$.  (It will be a multiple of the former
when $i<j$ and the latter when $i>j$).  This means that all determinants
in the sum on the right-hand side are zero except for the one in
which $i=j$. We get
$$
\left|\begin{matrix}
    h_{i-1} & h_{i} & h_{i+1} \\
    r_{i-1} & r_i & r_{i+1} \\
    r_{c-i+1} & r_{c-i} & r_{c-i-1}
  \end{matrix}\right|
= f_i \left|\begin{matrix}
    h_{\pi_i,i-1} & h_{\pi_i,i} & h_{\pi_i,i+1} \\
    r_{i-1} & r_i & r_{i+1} \\
    r_{c-i+1} & r_{c-i} & r_{c-i-1}
  \end{matrix}\right|
$$
and since $(h_{\pi_i,i-1},h_{\pi_i,i}, h_{\pi_i,i+1}) =
(r_{i-1},r_i,r_i r_{c-i-1}/r_{c-i})$ we get, by Laplace expansion
along the third column, that the determinant on the right-hand side
equals
$$
(\frac{r_ir_{c-i-1}}{r_{c-i}}-r_{i+1}) \left|\begin{matrix}
    r_{i-1} & r_i \\
    r_{c-i+1} & r_{c-i}
  \end{matrix}\right|
= \frac{1}{r_{c-i}}\left|\begin{matrix}
    r_{i} & r_{i+1} \\
    r_{c-i-1} & r_{c-i}
  \end{matrix}\right|
\left|\begin{matrix}
    r_{i-1} & r_i \\
    r_{c-i+1} & r_{c-i}
  \end{matrix}\right|.
$$
\end{proof}

\begin{thm} \label{NormHVec} Let $h(t) = \sum_{i \in \mathbb{Z}}
  h_it^i$ be a polynomial. 
Then $h(t)$ is a rational multiple of the $h$-vector of a level
$R$-module of codimension $\pdim$, generated in degree zero,
 whose Betti diagram is a non-negative linear
  combination of pure diagrams, if and only if
$$
\left|\begin{matrix}
    h_{i-1} & h_{i} & h_{i+1} \\
    r_{i-1} & r_i & r_{i+1} \\
    r_{c-i+1} & r_{c-i} & r_{c-i-1}
  \end{matrix}\right| \geq 0
$$
for all integers $i$, where $r_i = \binom{\pdim-1+i}{\pdim-1}$.
\end{thm}

\begin{proof}
  Assume first that the determinants of the assertion are non-negative
  for all integers $i$. Then, using Lemma~\ref{ECompLemma}, 
we get that $h(t) = \sum_{i=0}^c f_i h_{\pi_i}(t)$.
It follows from a straightforward calculation and the fact that
$r_i$ is a increasing function of $i$, that
$$
\left|\begin{matrix}
      r_{i-1} & r_i \\
      r_{c-i+1} & r_{c-i}
    \end{matrix}\right|
  \left|\begin{matrix}
      r_{i} & r_{i+1} \\
      r_{c-i-1} & r_{c-i}
    \end{matrix}\right| \geq 0
$$
and hence that the numbers $f_i$ and the determinants of the
assertion have the same sign. Hence, $f_i$ is positive for $i=0,1,\dots,c$. 
Let $M_i$ be an
  extremely compressed level module with Betti diagram $m_i\pi_i$ for
  some integer $m_i$, whose existence is guaranteed by
  Proposition~\ref{CompExistProp} and choose an integer $m$ such that
  the numbers $q_i = \frac{f_im}{m_i}$ are integers for $i =
  0,1,\dots,c$.  Then the module
$$
M = \bigoplus_{i=0}^c M_i^{q_i},
$$
where $M_i^{q_i}$ is the direct sum of $M_i$ with itself $q_i$ number
times, is level, its $h$-vector is equal to $mh(t)$ and its Betti
diagram, $\beta(M) = m\sum_{i=0}^c f_i\pi_i$, is a non-negative linear
combination of pure diagrams.

Assume now that $m$ is an integer such that $mh(t)$ is the $h$-vector
of a level $R$-module, $M$, whose Betti diagram is a non-negative
linear combination of pure diagrams. Then by
Proposition~\ref{MaxCompressedProp} we have
$$
\beta(M) = \sum_i f_i \pure(0,j+1,j+2,\dots,j+p-1,c+p) - \sum_{i,j}
b_{i,j}C^{i,j}
$$
for some non-negative rational numbers $f_i$. Since then $h_M(t) =
\sum_i f_i h_{\pi_i}(t)$, we get by Lemma~\ref{ECompLemma} that each
number $f_i$ have the same sign as the determinant in the proposition
which shows that these are non-negative.
\end{proof}

\begin{rmk} \label{LevelCodim2Rmk}
  When the codimension is two we have that $r_i = i+1$ and hence
$$
\left|\begin{matrix}
    h_{i-1} & h_{i} & h_{i+1} \\
    r_{i-1} & r_i & r_{i+1} \\
    r_{c-i+1} & r_{c-i} & r_{c-i-1}
  \end{matrix}\right|
= \left|\begin{matrix}
    h_{i-1} & h_{i} & h_{i+1} \\
    i & i+1 & i+2 \\
    c-i+2 & c-i+1 & c-i
  \end{matrix}\right|
= (c+2)\left|\begin{matrix}
    h_{i-1} & h_{i} & h_{i+1} \\
    0 & 1 & 2 \\
    1 & 1 & 1
  \end{matrix}\right|,
$$
where the second equality follows from the determinant on its
left-hand side by row operations on the second and third row.  Laplace
expansion along the first row of the rightmost determinant above gives
$$
h_{i-1}\left|\begin{matrix} 1 & 1 \\ 1 & 2 \end{matrix}\right| -
h_{i}\left|\begin{matrix} 1 & 1 \\ 0 & 2 \end{matrix}\right|+
h_{i+1}\left|\begin{matrix} 1 & 1 \\ 0 & 1 \end{matrix}\right| =
h_{i-1}-2h_i+h_{i+1}.
$$
We see that if the codimension is two then the condition given in
Theorem~\ref{NormHVec} is equal to $h_{i-1}-2h_i+h_{i+1} \geq 0$
for all $i$.  This is in fact precisely what characterizes the
$h$-vectors of level modules of codimension two as was shown in
\cite{Soderberg2}.
\end{rmk}

\section{Cancellable $h$-vectors}
Geramita and Lorenzini \cite{Geramita-Lorenzini} introduced the notion
of a cancellable $h$-vector of an algebra and this notion was
generalized to $h$-vectors of modules by the author in
\cite{Soderberg1} in order to use dualization to study level algebras.
We will now describe the cancellable $h$-vectors up to multiplication
by a rational number.

We begin by explaining the notion of a cancellable $h$-vector.
Consider the polynomial $h(t) = h_0 + h_1t + \dots + h_ct^c$ and
assume that it is the $h$-vector of a Cohen-Macaulay $R$-module $M$.
Then we have
\begin{equation} \label{CancelEq} \beta(M) = \MaxLex(h(t)) -
  \sum_{i,j} b_{i,j}C^{i,j}
\end{equation}
for some non-negative integers $b_{i,j}$. If $M$ is level, of
codimension $p$ and generated in degree zero, then its Betti diagram has
only one non-zero entry in column $p$, that is, $\beta(M)_{p,j} = 0$
whenever $j \neq c+p$.  We get, by considering the entries of position
$(p,j)$ and $(p-1,j)$ of (\ref{CancelEq}),
$$
\beta(M)_{p,j} = \MaxLex_{p,j}(h(t)) - b_{p-1,j}
$$ 
and
$$
\beta(M)_{p-1,j} = \MaxLex_{p-1,j}(h(t)) - b_{p-1,j} -b_{p-2,j}.
$$
Since $\beta(M)_{p,j} = 0$, whenever $j \neq c+p$, and
$\beta(M)_{p-1,j}$ and $b_{p-2,j}$ both are non-negative, this
means that
$$
\MaxLex_{p-1,j}(h(t)) - \MaxLex_{p,j}(h(t)) \geq 0
$$
for all $j \neq c+p$. Any $h$-vector that satisfies the above
inequality is said to be \emph{cancellable}.

In \cite{Soderberg1} the Eliahou-Kervaire resolution
\cite{Eliahou-Kervaire} was used to compute the difference
$\MaxLex_{p-1,j}(h(t)) - \MaxLex_{p,j}(h(t))$ in terms of the
$h$-vector $h(t)$. Due to the combinatorial nature of lexicographic
submodules the expression so obtained was rather cumbersome (see
Remark~\ref{OldCancelRmk}).  We will now see that the description of
cancellable $h$-vectors up to multiplication by a rational number, is
much simpler, and also much simpler to obtain.

\begin{prop} \label{BettiDiffProp} Let $h(t) = \sum_{i \in \mathbb{Z}}
  h_it^i$ be a polynomial of degree $c$ and let $D = \sum_{j=0}^c a_j
  \pure(0,j+1,j+2,\dots,j+p)$, where $a_j = h_j/r_j -
  h_{j+1}/r_{j+1}$. Then
$$
D_{p-1,p+j} - D_{p,p+j} = \frac{1}{r_{j+2}} \left| \begin{matrix}
    h_{j} & h_{j+1} & h_{j+2} \\
    r_{j} & r_{j+1} & r_{j+2} \\
    p & 1 & 0
  \end{matrix} \right|
$$
for all $j$ and $h(t)$ is a rational multiple of a cancellable
$h$-vector, of codimension $p$, if and only if $D_{p-1,p+j} - D_{p,p+j} \geq 0$, for $j =
0,1,\dots,c-1$, and $a_j \geq 0$, for $j=0,1,\dots,c$.
\end{prop}

\begin{proof}
  To compute the entries of the diagram $D$ we need to compute the
  entries of the pure diagrams $\pi_j = \pure(0,j+1,j+2,\dots,j+p)$
  for each integer $j$.  Herzog and K\"{u}hl calculated the Betti
  numbers of a pure resolution of any type \cite{Herzog-Kuhl}, and
  their calculation can be applied directly to pure diagrams (see
  \cite{Boij-Soderberg}). Their result states that
$$
\pure(d)_{k,d_k} = (-1)^k \prod_{i\neq k} \frac{d_i}{d_i-d_k}
$$
for any type $d = (0,d_1,d_2,\dots,d_p)$.  When $d$ is the type of
$\pi_j$ we have that $d_i = j+i$ for each $i = 1,2,\dots,p$.  The
formula of Herzog and K\"{u}hl now yields
\[
  (\pi_j)_{p,p+j} = (-1)^p \prod_{i\neq p} \frac{d_i}{d_i-d_p} =
  \prod_{i=1}^{p-1} \frac{j+i}{p-i} = \frac{(j+p-1)!}{j!(p-1)!} =
  \binom{p-1+j}{p-1} = r_j
\]
and
\[\begin{split}
  (\pi_j)_{p-1,p+j-1}& = (-1)^p \prod_{i\neq p-1} \frac{j+i}{i-p+1} =
  (j+p)\prod_{i=1}^{p-2} \frac{j+i}{p-i-1} \\
  & =(j+p)\frac{(p+j-2)!}{j!(p-2)!}  =
  \frac{(j+p)(p-1)r_j}{p+j-1} 
  = pr_j-r_{j-1}.
\end{split}\]
The entries of $D$ that we need are $D_{p,p+j}$ and $D_{p-1,p+j}$ and
these entries are given by $D_{p,p+j} = a_j(\pi_j)_{p,p+j}$ and
$D_{p-1,p+j} = a_{j+1}(\pi_{j+1})_{p-1,p+j}$. By the above calculation
we get
$$
D_{p-1,p+j} - D_{p,p+j} =  a_{j+1}(pr_{j+1} - r_{j}) - a_jr_j \\
$$
which, using the assumption that $a_j = h_j/r_j - h_{j+1}/r_{j+1}$, in
turn gives
\begin{equation*}
  \begin{split} 
    D_{p-1,p+j} - D_{p,p+j} =& \left(\frac{h_{j+1}}{r_{j+1}} -
      \frac{h_{j+2}}{r_{j+2}}\right)(-r_{j} + pr_{j+1}) -
    \left(\frac{h_j}{r_j} - \frac{h_{j+1}}{r_{j+1}}\right) r_j \\
    =& -h_j + ph_{j+1} - \frac{pr_{j+1} - r_j}{r_{j+1}}h_{j+2}.
  \end{split}
\end{equation*}
To finish the computation of $D_{p-1,p+j} - D_{p,p+j}$, note that
$$
\frac{1}{r_{j+2}} \left| \begin{matrix}
    h_{j} & h_{j+1} & h_{j+2} \\
    r_{j} & r_{j+1} & r_{j+2} \\
    p & 1 & 0
  \end{matrix} \right| = -h_j + ph_{j+1} - \frac{pr_{j+1} -
  r_j}{r_{j+1}}h_{j+2}.
$$

By Proposition~\ref{SingleMaxBettiProp}, $h(t)$ is an $h$-vector if
and only if $a_j \geq 0$, for $j=0,1,\dots,c$, and furthermore, then it exists
 an integer $s$ such that $\MaxLex(sh(t)) = sD$.  Assume that
$D_{p-1,p+j} - D_{p,p+j} \geq 0$, for $j=0,1,\dots,c-1$ and note that
since $D_{p,p+j} = 0$ when $j < 0$ of $j>c$ this means that
$D_{p-1,p+j} - D_{p,p+j} \geq 0$, for all $j \neq c$.  Then, by
definition, $sh(t)$ is cancellable and we see that $h(t)$ is a
rational multiple of the cancellable $h$-vector $sh(t)$.

If $h(t)$ is cancellable then, as we will see, any integer multiple of
$h(t)$ is cancellable and in particular $sh(t)$ which shows that
$D_{p-1,p+j} - D_{p,p+j} \geq 0$ for all $j \neq c$.  Let $m$ be an
integer and consider the the polynomial $mh(t)$.  Then
$m\MaxLex((h(t))$ is the Betti diagram of $(F/L)^m$ where $L$ is the
lexicographic submodule such that $F/L$ have $h$-vector $h(t)$, and
since $(F/L)^m$ have $h$-vector $mh(t)$ we get
$$
m\MaxLex(h(t)) = \MaxLex(mh(t)) - \sum_{i,j}b_{i,j}C^{i,j}.
$$
for some non-negative integers $b_{i,j}$.  By considering the entries
in position $(p-1,j)$ and $(p,j)$ of the equation above, we get
$$
m\MaxLex_{p,j}(h(t)) = \MaxLex_{p,j}(mh(t)) - b_{p-1,j}
$$ 
and
$$
m\MaxLex_{p-1,j}(h(t)) = \MaxLex_{p,j}(mh(t)) - b_{p-1,j} -b_{p-2,j}
$$
from which it follows that
$$
m\MaxLex_{p-1,j}(h(t)) - m\MaxLex_{p,j}(h(t)) = \MaxLex_{p,j}(mh(t)) -
\MaxLex_{p,j}(mh(t)) - b_{p-2,j}.
$$
Since $b_{p-2,j}$ is non-negative this means that
$$
\MaxLex_{p-1,j}(mh(t)) - \MaxLex_{p,j}(mh(t)) \geq
m\left(\MaxLex_{p-1,j}(h(t)) - \MaxLex_{p,j}(h(t))\right)
$$
for all $j$, which shows that if $h(t)$ is cancellable then $mh(t)$ is
cancellable as well.
\end{proof}

The following proposition uses that the dual of a level module is level
and hence that they both have cancellable $h$-vectors which then
satisfy the condition of \ref{BettiDiffProp}.  For any polynomial $h_0
+ h_1t + \dots + h_ct^c$ we call the polynomial $h_c + h_{c-1}t +
\dots + h_0t^c$ its \emph{reverse}, and we note that if $M$ is a
Cohen-Macaulay module with $h$-vector $h(t)$ of degree $c$, then the reverse of
$h(t)$ is the $h$-vector of $M^{\vee}(-c)$.

\begin{thm} \label{DualCancelProp}Let $h(t) = \sum_{i \in \mathbb{Z}}
  h_it^i$ be a polynomial of degree $c$. Then $h$ is a
  rational multiple of a cancellable $h$-vector whose reverse also is
  cancellable if and only if
$$
\left| \begin{matrix}
    h_{i-1} & h_{i} & h_{i+1} \\
    r_{i-1} & r_{i} & r_{i+1} \\
    p & 1 & 0
  \end{matrix} \right| \geq 0, \quad \left|
  \begin{matrix}
    h_{i-1} & h_{i} & h_{i+1} \\ 0 & 1 & p \\
    r_{c-i+1} & r_{c-i} & r_{c-i-1}
  \end{matrix} \right| \geq 0
$$
for all $i=1,2, \dots,c$, and
$$
\frac{h_i}{r_i}-\frac{h_{i+1}}{r_{i+1}} \geq 0, \quad 
\frac{h_{i+1}}{r_{c-i-1}}-\frac{h_{i}}{r_{c-i}} \geq 0
$$
for all $i=0,1, \dots,c$, where $r_i = \binom{p-1+i}{p-1}$, and as
a consequence all $h$-vectors of level modules, generated in degree
zero, satisfy these two conditions.
\end{thm}

\begin{proof}
  If $h$ and its reverse are cancellable then conditions of
  the proposition are satisfied by applying
  Proposition~\ref{BettiDiffProp} to $h$ and its reverse.

  Assume now that these conditions are satisfied. Then, by
  Proposition~\ref{BettiDiffProp}, there are integers $q$ and
  $q'$ such that $qh(t)$ and $q'h'(t)$, where $h'(t)$ is the reverse
  of $h(t)$, both are cancellable. Since, as shown in the last part of
  the proof of Proposition~\ref{BettiDiffProp}, this means that
  $qq'h(t)$ and $qq'h'(t)$ both are cancellable, we see that $h(t)$ is
  a rational multiple of a cancellable $h$-vector whose reverse is also
   cancellable.

  If $h(t)$ is the $h$-vector of a level module $M$, then $h$ is
  cancellable. Since $M^{\vee}(-c)$ is level as well, and
  furthermore, generated in degree zero and have the reverse of $h(t)$
  as its $h$-vector, we see that the the reverse of $h(t)$ is
  cancellable.
\end{proof}

\begin{rmk} \label{OldCancelRmk}
  By using the Eliahou-Kervaire resolution \cite{Eliahou-Kervaire} the
  following expression is obtained in \cite[Proposition
  17]{Soderberg1}
  \begin{multline} \label{OldDiff}
    \MaxLex_{p-1,p+j}(h(t)) - \MaxLex_{p,p+j}(h(t)) = \\
    p\left(h_{j+1}-q r_{i+1}-\MacRep{s}{j+2}{-1}{-1}\right) -h_{j}+q
    r_{j}+\MacRep{s}{j+2}{-2}{-2}
  \end{multline}
  where $q$ is the quotient and $s$ the remainder when $h_{j+2}$ is
  divided by $r_{j+2}$ and the expressions $\MacRep{s}{j+2}{-1}{-1}$
  and $\MacRep{s}{j+2}{-2}{-2}$ are manipulations with the $(j+2)$-th
  Macaulay representation of $s$.  We will now rearrange the
  right-hand side of this equality and see that if we choose the right
  integer multiple of $h(t)$ we get the result of
  Proposition~\ref{BettiDiffProp}.
  The right-hand side of (\ref{OldDiff}) equals, after some
  rearranging,
$$
-h_j + ph_{j+1} - q(pr_{j+1} - r_j)
+\MacRep{s}{j+2}{-2}{-2}-p\MacRep{s}{j+2}{-1}{-1}
$$
and using $q = (h_{j+2}-s)/r_{j+2}$, we see that this equals
\begin{equation} \label{oldEq} -h_j + ph_{j+1} - \frac{pr_{j+1} -
    r_j}{r_{j+2}}h_{j+2} -\left(
    -\MacRep{s}{j+2}{-2}{-2}+p\MacRep{s}{j+2}{-1}{-1}-\frac{pr_{j+1} -
      r_j}{r_{j+2}}s\right).
\end{equation}
Choose an integer $m$ such that $mh_{j+2}$ is divisible by $r_{j+2}$
for each $j$.  Then the number $s$ of (\ref{oldEq}), which is the
remainder when $mh_{j+2}$ is divided by $r_{j+2}$, is $s=0$, and hence
\[
\MaxLex_{p-1,p+j}(mh(t)) - \MaxLex_{p,p+j}(mh(t)) = m\left(-h_j +
  ph_{j+1} - \frac{pr_{j+1} - r_i}{r_{j+2}}h_{j+2} \right).
\]
The expression above is $m$ times the difference
$D_{p-1,p+j}-D_{p,p+j}$ of Proposition~\ref{BettiDiffProp}, so we get
in this way the same description of cancellable $h$-vectors up to
multiplication by a rational numbers as the one in given in that
proposition.
\end{rmk}

We will now see how the conditions of Theorem~\ref{NormHVec} and
Theorem~\ref{DualCancelProp} behave in some cases.  In what
follows we restrict the notion of a cancellable $h$-vector to those
which also have a cancellable reverse. If $h(t)$ is the $h$-vector of
a level module we say that $h(t)$ is level.

Geramita et al. \cite{Geramita-Harima-Migliore-Shin} found all
$h$-vectors of artinian level algebras of codimension three which have 
either socle degree less than six or socle degree six
and type two (the ones of type one are Gorenstein so they are known
for any socle degree by Stanley's result \cite{Stanley}) If $M$ is level then both
$h(t)$ and its reverse are $h$-vectors of modules generated in a
single degree, so they both satisfy Macaulay's condition for modules
generated in a single degree. Considering the set of $h$-vectors which
satisfy, and also have a reverse that satisfy, the condition of
Macaulay, we will now see, in one of the cases covered by Geramita et
al. how many are level, cancellable and are a rational multiple of a
cancellable $h$-vector.  We will also see how many of them satisfy the
condition of Theorem~\ref{NormHVec}.  Remember that all
polynomials that satisfy the condition in Theorem~\ref{NormHVec}
are rational multiples of level $h$-vectors, and we conjecture the
converse to hold, that is, that all $h$-vectors that are rational
multiples of level $h$-vectors satisfy the condition of
Theorem~\ref{NormHVec}.

\begin{example}
  There are 148 polynomials on the form $h(t) = 1 + h_1t + h_2t^2 +
  \dots + h_5t^5 + 2t^6$, with non-negative integer coefficient, which
  satisfy, and also have a reverse that satisfy, Macaulay's condition
  for modules generated in a single degree.  Among these are, as
  mentioned, the level ones and there are 58 of them.  Of the 148, 67
  satisfy the condition given in Theorem~\ref{NormHVec}, that is,
  they are non-negative linear combinations of $h$-vectors of
  extremely compressed level modules of some fixed socle degree, and
  all the level ones are among these, which if
  Conjecture~\ref{MainConj} is true always will be the case.  There
  are nine $h$-vectors, among these 67, which are not level and these
  are
  \[\begin{split}
    1+3t+4t^2+5t^3+6t^4+4t^5+2t^6,&\quad
    1+3t+4t^2+5t^3+6t^4+6t^5+2t^6, \\
    1+3t+5t^2+5t^3+4t^4+3t^5+2t^6,&\quad
    1+3t+5t^2+6t^3+7t^4+4t^5+2t^6, \\
    1+3t+5t^2+7t^3+7t^4+4t^5+2t^6,&\quad
    1+3t+5t^2+7t^3+9t^4+5t^5+2t^6, \\
    1+3t+6t^2+5t^3+4t^4+3t^5+2t^6,&\quad
    1+3t+6t^2+6t^3+5t^4+4t^5+2t^6, \\
    1+3t+6t^2+10t^3+7t^4+5t^5+2t^6.
  \end{split}\]
  However, since these nine $h$-vectors satisfy the condition in
  Theorem~\ref{NormHVec} they are rational multiples of level
  $h$-vectors.  Of the 148 that satisfy Macaulay's condition for
  modules generated in a single degree, 71 are cancellable and 116
  become cancellable after multiplication by some rational number.

  It turns out that the 67 $h$-vectors satisfying the condition of
  Theorem~\ref{NormHVec}, they are then rational multiples of
  level $h$-vectors and hence of cancellable $h$-vectors, are all
  cancellable and the cancellable ones that does not satisfy the
  condition of Theorem~\ref{NormHVec} are
  \[\begin{split}
    1+3t+5t^2+7t^3+6t^4+6t^5+2t^6, &\quad
    1+3t+6t^2+7t^3+6t^4+6t^5+2t^6, \\
    1+3t+6t^2+9t^3+7t^4+6t^5+2t^6, &\quad
    1+3t+6t^2+10t^3+7t^4+6t^5+2t^6. \\
  \end{split}\]
  Considering all the the socle degrees and types covered by Geramita
  et al. there are only three more examples of cancellable $h$-vectors
  that does not satisfy the condition of Theorem~\ref{NormHVec},
  and these are
  \[\begin{split}
    1+3t+5t^2+6t^3+4t^4+3t^5, &\quad
    1+3t+6t^2+6t^3+4t^4+3t^5, \\
    1+3t+6t^2+10t^3+7t^4+6t^5.
  \end{split}\]

  Again considering all the cases covered by Geramita et al., there
  are eight $h$-vectors that satisfy the condition of
  Theorem~\ref{NormHVec} but which are not cancellable.  They are
  \[\begin{split}
    1+3t+6t^2+8t^3+8t^4+9t^5, &\quad
    1+3t+6t^2+9t^3+9t^4+10t^5, \\
    1+3t+6t^2+8t^3+9t^4+11t^5, &\quad
    1+3t+6t^2+9t^3+10t^4+12t^5, \\
    1+3t+6t^2+10t^3+11t^4+13t^5, &\quad
    1+3t+6t^2+10t^3+12t^4+15t^5, \\
    1+3t+6t^2+6t^3+7t^4, &\quad 1+3t+6t^2+7t^3+9t^4.
  \end{split}\]
\end{example}

We now take a closer look at one of the $h$-vector above which
satisfies the condition of Theorem~\ref{NormHVec}, but which is not
cancellable.

\begin{example}
  The polynomial $h(t) = 1+3t+6t^2+7t^3+9t^4$ satisfies the condition
  of Theorem~\ref{NormHVec} and hence is a rational multiple of
of some level $h$-vector, but it is not cancellable. In fact, the lexicographic
  ideal $L$ of $R=k[x,y,z]$ such that $R/L$ have the polynomial $h(t)$
  as $h$-vector is
$$
L = (z^3, yz^2, xz^2, y^4z, xy^3z, x^2y^2z, x^3yz, x^4z, y^5, xy^4,
x^2y^3, x^3y^2, x^4y, x^5)
$$
and hence $\MaxLex(h(t)) = \beta(R/L)$ and this diagram is in fact
given by
$$
\MaxLex(h(t)) = \begin{pmatrix} 1&-&-&- \\ - &-&-&- \\ -&3&3&1 \\ -
  &-&-&- \\ -&11&20&9 \end{pmatrix},
$$
where we have used the convention of writing the entry
$\MaxLex_{i,j}(h(t))$ in column $i$ and row $j-i$.  We see that $h(t)$
is not cancellable, and hence not level, since there is no way to
cancel the Betti number $\beta_{3,5}=1$ in the last column of
$\MaxLex(h(t))$.

The diagram from Proposition~\ref{SingleMaxBettiProp} that bounds from
above the Betti numbers of modules with $h$-vector $h(t)$ is
\begin{multline*}
D = \frac{3}{10}\beta(R/\mathfrak{m}^3) + \frac{1}{10}\beta(R/\mathfrak{m}^4) +
  \frac{3}{5}\beta(R/\mathfrak{m}^5)=
  \begin{pmatrix}
    1&-&-&- \\ -&-&-&-\\ -& 3 & 9/2 & 9/5 \\ -& 3/2& 12/5& 1\\ -&
    63/5& 21& 9 \end{pmatrix},
\end{multline*}
and since $\beta_{3,5}=9/5 < \beta_{2,5} = 12/5$ and $\beta_{3,6} = 1
< \beta_{2,6} = 21$ we see that $h(t)$ is a rational multiple of a
cancellable $h$-vector. Indeed, 
$\MaxLex(10h(t)) = 10 D$, so $10h(t)$ is cancellable.

In general it is hard to say if a cancellation of a Betti diagram can
be realized as the Betti diagram of some module, but in this case,
since $h(t)$ is a rational multiple of a level $h$-vector, we know
that there is a cancellation of a rational multiple of $D$ such that
$\beta_{3,5} = \beta_{3,6} = 0$.  We also know that $h(t)$ is a
non-negative linear combination of $h$-vectors of extremely compressed
level modules of some fixed socle degree, since it satisfies the
condition of Theorem~\ref{NormHVec}. In fact,
\begin{multline*}
  h(t) = 1+3t+6t^2+7t^3+9t^4  = \\
  \frac{3}{7}(1+3t+6t^2+3t^3+t^4) + \frac{4}{7}(1+3t+6t^2+10t^3+15t^4)
\end{multline*}
where $1+3t+6t^2+3t^3+t^4$ and $1+3t+6t^2+10t^3+15t^4$ are the
$h$-vectors of some extremely compressed level modules, $M_1$ and
$M_2$, with Betti diagrams $\pure(0,3,4,7)$ and $\pure(0,5,6,7)$,
respectively. The existence of these extremely compressed modules are
guaranteed by Proposition~\ref{CompExistProp}.  The corresponding
diagram of this linear combination is
\begin{multline*}
  E = \frac{3}{7} \pure(0, 3, 4, 7) + \frac{4}{7} \pure(0, 5, 6, 7) = \\
  \frac{3}{7} \begin{pmatrix} 1&-&-&- \\ -&-&-&- \\ -&7&7&- \\ -&-&-&-
    \\ -& -& -& 1 \end{pmatrix} + \frac{4}{7} \begin{pmatrix} 1&-&-&-
    \\ -&-&-&- \\ -&-&-&- \\ -&-&-&- \\ -& 21& 35& 15 \end{pmatrix} =
  \begin{pmatrix} 1&-&-&- \\ -&-&-&- \\ -&3&3&- \\ -&-&-&- \\ -& 12&
    20& 9 \end{pmatrix}.
\end{multline*}
We see that $7 E$, which by Proposition~\ref{SingleMaxBettiProp} is
obtained from the diagram $7D$ by a sequence of consecutive
cancellations, is the Betti diagram of a level module with $h$-vector
$7 h(t)$.
\end{example}

\section{The weak Lefschetz property}
In \cite{Harima-Migliore-Nagel-Watanabe} all the Hilbert functions of
artinian algebras with the weak or strong Lefschtez property are
determined and a sharp upper bound on the graded Betti numbers of
these algebras is given. We give analogous results in the case of
level modules with the weak Lefschetz property, but only up to
multiplication by a rational number.  Furthermore, we show that the
upper bound on the graded Betti numbers, in our case, is a
non-negative linear combination of pure diagrams and that the graded
Betti numbers are obtained from these by a sequence of consecutive
cancellations.  In Section \ref{CodimThree} we will see that this
upper bound is enough to prove Conjecture~\ref{MainConj} for level
modules of codimension three with the weak Lefschetz property, which
in turn gives the Multiplicity conjecture for these modules.

An artinian $R$-module $M$ has the weak Lefschetz property if there is
a linear form $\ell$ such that $M_{i-1} \xrightarrow{\times \ell}
M_{i}$ is either injective or surjective for each $i$, and such a form
is called a Lefschetz element. Let $M$ be a graded $R$-module and
$\ell$ a linear form of $R$, not necessarily a Lefschetz element.  Then
we have an exact sequence
$$
0 \to P \to M(-1) \xrightarrow{\times \ell} M \to Q \to 0,
$$
where the map from $M(-1)$ to $M$ is multiplication by $\ell$, and $P$
and $Q$ are the kernel and cokernel, respectively. Then $P$ and $Q$
are modules over $S = R/\ell$ and hence the codimensions of $P$ and $Q$
are the codimension of $M$ minus one. In \cite{Migliore-Nagel}, Migliore and
Nagel shows that the exact sequence above gives the following exact
sequences of graded $R$-modules
\begin{multline} \label{LongTorEq}
  \dots \to \Tor^S_{i-1}(P,k) \to \Tor^R_i(M,k) \to \Tor^S_i(Q,k) \to \\
  \dots \to \Tor^S_{0}(P,k) \to \Tor^R_1(M,k) \to \Tor^S_1(Q,k) \to 0
\end{multline}
and
\begin{equation} \label{TorBeginEq} 0 \to \Tor^S_{0}(M,k) \to
  \Tor^R_0(Q,k) \to 0.
\end{equation}
Note that the long exact sequence above yields
\begin{equation} \label{TorEndEq} 0 \to \Tor^S_{p-1}(P,k) \to
  \Tor^R_p(M,k) \to 0
\end{equation}
These exact sequences implies, as observed by Migliore and Nagel, that
\[ 
\beta^R_{i,j}(M) \leq
  \beta^S_{i,j}(Q) +\beta^S_{i-1,j}(P)
\]
for all $i$ and $j$, with equality when $\beta^S_{i-1,j}(P) =
\beta^S_{i-2,j}(P) = 0$ or $\beta^S_{i+1,j}(Q) = \beta^S_{i,j}(Q) =
0$.  We now make a further observation. Denote by $D$ the diagram
defined by
$$
D_{i,j} = \beta^S_{i,j}(Q) +\beta^S_{i-1,j}(P).
$$
for all $i$ and $j$.
The exact sequence (\ref{LongTorEq}) then yields
$$
\sum_{k=1}^i (-1)^{k+1} (D_{k,j} - \beta_{k,j}(M)) \geq 0,
$$
for all $i$, since $\beta_{i,j}(M) = \dim_k \Tor_i(M,k)_j$ and
$\dim_k$ is additive along exact sequences. 
This implies not only that $D \leq \beta(M)$, entry by
entry, but also that $\beta(M)$ is obtained from $D$ by a sequence of
consecutive cancellations. In fact, let
$$
b_{i,j} = \sum_{k=1}^{i-1} (-1)^{k+1} (D_{i,j} - \beta_{i,j}(M)).
$$
Then
$$
\beta(M) = \beta(D) - \sum_{i,j}b_{i,j}C^{i,j}.
$$

\begin{lemma} \label{LefLemma} Let $M$ be a level artinian $R$-module
  of codimension $p$, generated in degree zero, with a Lefschetz
  element $\ell$ and let
$$
0 \to P \to M(-1) \xrightarrow{\times \ell} M \to Q \to 0,
$$
be the exact sequence corresponding to the multiplication map given by
$\ell$. Let $h(t) = \sum_{i \in \mathbb{Z}}h_it^i$ be the $h$-vector
of $M$ and assume that it is a polynomial of degree $c$.  Furthermore,
let $u$ be the smallest integer such that $h_{u} \geq h_{u+1}$.  Then
the dual of $P$, $P^{\vee}$, is generate in a single degree and that
degree is $-c-1$, and hence $P^{\vee}(-c-1)$ is generated in degree
zero. Furthermore, the $h$-vector of $Q$ is
$$
h_Q(t) = \sum_{i=0}^u (h_i-h_{i-1})t^i,
$$ 
the $h$-vector of $P$ is,
$$
h_P(t) = \sum_{i = u+1}^{c+1} (h_{i-1} - h_{i})t^i.
$$
the $h$-vector of $P' = P^{\vee}(-c-1)$ is,
$$
h_{P'}(t) = \sum_{i=0}^{c-u} (h_{c-i}-h_{c-i+1})t^i
$$
and for the multiplicities it holds that
$$
e(Q) = e(P) = e(P^{\vee}(-c-1)).
$$
\end{lemma}

\begin{proof}
  The socle degree of $M$ is $c$ since the degree of $h(t)$ is $c$.
  This means that $\Tor^R_p(M,k)$ is concentrated in degree $c+p$.  It
  follows from (\ref{TorEndEq}) that $\Tor^S_{p-1}(P,k) \cong
  \Tor^R_p(M,k)$ and hence $\Tor^S_{p-1}(P,k)$ is also concentrated in
  degree $c+p$, which shows that $\beta_{p-1,j}(P) = 0$ for all $j
  \neq c+p$. The graded Betti numbers of the dual of $P$ are obtained
  from the ones of $P$ by $\beta_{i,j}(P) =
  \beta_{p-1-i,p-1-j}(P^{\vee})$, since $P$ is of codimension $p-1$,
  and this shows that $\beta(P^{\vee})_{0,j} = 0$ for all $j \neq
  -c-1$. Hence, $P^{\vee}$ is generated in degree $-c-1$.

  Since $u$ is the smallest integer such that $h_{u} \geq h_{u+1}$, we
  see that $u$ is also the smallest integer such that the map $M_{i}
  \to M_{i+1}$, given by multiplication by $\ell$, is surjective.
  Since $M$ is level, and hence generated in a single degree, we see
  that as soon as the map $M_{i+1} \to M_{i}$ is surjective it will
  continue to be in all higher degrees.  Hence $M_{i+1} \to M_{i}$ is
  surjective for all $i \geq u$ and injective for all $i < u$.  The
  exact sequence
$$
0 \to P_i \to M_{i-1} \xrightarrow{\times \ell} M_i \to Q_i \to 0,
$$
then implies that $Q_i = 0$ for all $i > u$ and $P_i = 0$ for all $i
\leq u$.  It also follows from the exact sequence that
\begin{equation}\label{HilbPQMEq}
  H(Q,i) - H(P,i) = H(M,i) - H(M,i-1) = h_{i} - h_{i-1},
\end{equation}
and using that $Q_i = 0$ for all $i > u$ and that $P_i = 0$ for all $i
\leq u$, we get
$$
h_Q(t) = \sum_{i \in \mathbb{Z}} H(Q,i)t^i = \sum_{i = 0}^{u} (h_{i} -
h_{i-1})t^i
$$
and
$$
h_P(t) = \sum_{i \in \mathbb{Z}} H(P,i)t^i = \sum_{i = u+1}^{c+1}
(h_{i-1} - h_{i})t^i.
$$
Now, since $H(P^{\vee}(-c-1),i) = H(P^{\vee},i-c-1) = H(P,c+1-i)$ we
get that the $h$-vector of $P' = P^{\vee}(-c-1)$ is
$$
h_{P'}(t) = \sum_{i \in \mathbb{Z}} H(P^{\vee}(-c-1),i)t^i = \sum_{i
  \in \mathbb{Z}} H(P,c+1-i)t^i = \sum_{i=0}^{c-u}
(h_{c-i}-h_{c+1-i})t^i.
$$

It remains to prove that $e(Q) = e(P) = e(P^{\vee}(i-c-1))$. From
(\ref{HilbPQMEq}) it follows that
$$
h_Q(t) - h_P(t) = (1-t)h_M(t)
$$
and if we put $t=1$ in this equation we get $h_Q(1) = h_P(1)$.  By
definition, $e(Q) = h_Q(1)$ and $e(P) = h_P(1)$ so we get that $e(Q) =
e(P)$. Since the set of coefficients of $h_P(t)$ and $h_P'(t)$ are the
same we see that $h_{P'}(1) = h_P(1)$ which shows that $e(P) = e(P') =
e(P^{\vee},i-c-1))$.
\end{proof}

We need the following lemma which is almost a restatement of
\cite[Lemma 4.1]{Boij-Soderberg}.

\begin{lemma} \label{StepUpLemma} Let $F$ and $G$ be non-negative
  linear combinations of pure diagrams of codimension $p-1$ such that
  $e(F) = e(Q)$ and the maximal shifts of $F$ are are smaller than the
  minimal shifts of $G$, that is, $\up(F) < \down(G)$. Then the
  diagram $E$ defined by
$$
E_{i,j} = F_{i,j} + G_{i-1,j}
$$
for all integers $i$ and $j$, is a non-negative linear combination of
pure diagrams of codimension $p$.
\end{lemma}

\begin{proof}
  This follows from \cite[Lemma 4.1]{Boij-Soderberg} which states that
  the diagram $D$ defined by
$$
D_{i,j} = \pure(d)_{i,j} + \frac{e(\pure(d))}{e(\pure(d'))}\pure(d')_{i-1,j}
$$
for all integers $i$ and $j$, is a non-negative linear combination of
pure diagrams of codimension $p$, if the types of the codimension
$p-1$ diagrams $\pure(d)$ and $\pure(d')$ satisfies $d < d'$.  Assume
that $F = \sum_{i=1}^s F_i$ and $G = \sum_{i=1}^t G_i$ where $F_i$ and
$G_i$ are pure diagrams.  Then the type of $F_1$ is smaller than the
type of $G_1$, by the assumption on the minimal and maximal shifts of
$F$ and $G$.  Assume furthermore that $e(F_1) \leq e(G_1)$, the case
where $e(F_1) > e(G_1)$ follows in the same way, then the diagram
$E_1$ defined by
$$
(E_1)_{i,j} = (F_1)_{i,j} + \frac{e(F_1)}{e(E_1)}(G_1)_{i,j}
$$
for all integers $i$ and $j$, is non-negative linear combination of
pure diagrams by \cite[Lemma 4.1]{Boij-Soderberg}.  The diagram $E' =
E - E_1$ is now given by the diagrams $F' = \sum_{i=2}^s E_i$ and
$$
G' = \left( 1-\frac{e(F_1)}{e(E_1)}\right)G_1 + \sum_{i=2}^t G_i
$$
in the same way as $E$ is given by $F$ and $G$ in the statement of the
lemma.  The total number of pure diagrams in $F'$ and $G'$ are
$s+t-1$, that is, one fewer than the total number in $F$ and $G$, and
furthermore $e(F') = e(G')$.  By induction we may assume that $E'$ is
a non-negative linear combination of pure diagrams and then this is
true for $E$ as well since $E = E' + E_1$.
\end{proof}

\begin{prop} \label{LefProp} Let $M$ be a level artinian $R$-module of
  codimension $p$, generated in degree zero, with the weak
  Lefschetz property and let the $h$-vector of $M$ be the degree $c$
  polynomial $h(t) = \sum_{i \in \mathbb{Z}}h_it^i$.
  Furthermore, let $S=R/\ell$, $s_i = \binom{p-2+i}{p-2}$, $u$ be the
  smallest integer such that $h_u \geq h_{u+1}$ and
$$
F = \sum_{i=0}^u f_i \beta^S(S/\mathfrak{m}^{i+1})
$$
where
$$
f_i = \frac{h_i-h_{i-1}}{s_i} - \frac{h_{i+1} - h_{i}}{s_{i+1}}
$$
for $i = 0,1,\dots,u-1$ and
$$
f_u = \frac{h_u-h_{u-1}}{s_u}
$$
and
$$
G = \sum_{i=0}^{c-u} g_i \beta^S(S/\mathfrak{m}^{i+1})
$$
where
$$
g_i = \frac{h_{c-i} - h_{c-i+1}}{s_i} - \frac{h_{c-i-1} -
  h_{c-i}}{s_{i+1}}
$$
for $i=0,1,\dots,c-u-1$ and
$$
g_{c-u} = \frac{h_{u} - h_{u+1}}{s_{c-u}}.
$$
Then $f_i \geq 0$, for $i=0,1,\dots,u$, and $g_i \geq 0$, for $i =
0,1,\dots,c-u$, and $\beta(M)$ is obtained from the diagram $E$,
defined by
$$
E_{i,j} = F_{i,j} + G_{p-i,p+c-j},
$$ 
for all integers $i$ and $j$, by a sequence of consecutive
cancellations and furthermore, $E$ is non-negative linear combination
of pure diagrams of codimension $p$.
\end{prop}

\begin{proof}
  Let $\ell$ be a Lefschetz element on $M$, and consider the exact
  sequence
$$
0 \to P \to M(-1) \xrightarrow{\times \ell} M \to Q \to 0,
$$
where $P$ and $Q$ are the kernel and cokernel, respectively, of the
map $M(-1) \to M$ given by multiplication by $\ell$.  Then $P$ and $Q$
are artinian modules over $S = R/\ell$ and hence of codimension $p-1$.
The $h$-vector of $Q$ is, by Lemma~\ref{LefLemma}, $\sum_{i = 0}^u
(h_{i} - h_{i-1})t^i$ so it follows by
Proposition~\ref{SingleMaxBettiProp} that
\begin{equation}\label{FCancelEq}
  \beta^S(Q) = F - \sum_{i,j}b^F_{i,j}C^{i,j}.
\end{equation}
for some non-negative rational numbers $b^F_{i,j}$.  In the same way
we get, since the $h$-vector of $P^{\vee}(-c-1)$ is $\sum_{i=0}^{c-u}
(h_{c-i}-h_{c-i+1})t^i$ and $P^{\vee}(-c-1)$, by Lemma~\ref{LefLemma},
is generated in degree zero, that
\begin{equation}\label{EDualEq}
  \beta^S(P^{\vee}(-c-1)) = G - \sum_{i,j}b^G_{i,j}C^{i,j}.
\end{equation}
for some non-negative rational numbers $b^G_{i,j}$.

Since $P$ is of codimension $p-1$ we have by a standard property of
dualization that
$$
\beta^S_{i,j}(P) = \beta^S_{p-1-i,p-1-j}(P^{\vee})
$$
and since shifting the degrees of $P$ shifts the degrees of its Betti
numbers by the same amount we get
$$
\beta^S_{i,j}(P) = \beta^S_{p-1-i,p+c-j}(P^{\vee}(-c-1)).
$$
So, if we denote by $G'$ the diagram defined by
\begin{equation}\label{GEq}
  G'_{i,j} = G_{p-1-i,p+c-j}
\end{equation}
for all integers $i$ and $j$, we get, by (\ref{EDualEq}), that
$\beta^S(P)$ is obtained from $G'$ by sequence of consecutive
cancellations.  Let $D$ and $E$ be the diagrams defined by
$$
D_{i,j} = \beta^S_{i,j}(Q) + \beta^S_{i-1,j}(P)
$$ 
and
$$
E_{i,j} = F_{i,j} + G'_{i-1,j} = F_{i,j} + G_{p-i,c+p-j},
$$ 
for all integers $i$ and $j$.  Then, since $\beta^S(Q)$ and
$\beta^S(P)$ are obtained from $F$ and $G'$, respectively, by a
sequence of consecutive cancellations we see that $D$ is obtained from
$E$ in this way as well. As noted in the beginning of this section,
the long exact sequence (\ref{LongTorEq}) implies that $\beta(M)$ is
obtained from the diagram $D$ by a sequence of consecutive
cancellations, and hence $\beta(M)$ is be obtained from $E$ by a
sequence of consecutive cancellations.

It remains to prove that $E$ is a non-negative linear combination of
pure diagrams of codimension $p$.  We will prove this by applying
Lemma~\ref{StepUpLemma} to the diagrams $F$ and $G'$, and thus we have
to show that $F$ and $G$ are non-negative linear combinations of pure
diagrams and that their minimal and maximal shifts satisfies $\up(F) <
\down(G)$.

The Betti diagram of $S/\mathfrak{m}^{i+1}$, for any integer $i$, 
is the pure diagram $\pure(0,i+1,i+2,\dots,i+p-1)$ of codimension $p-1$, 
and since the rational numbers $f_i$ and $g_i$ are all non-negative we see
that the diagrams $F$ and $G$ are non-negative linear combinations
of pure diagrams of codimension $p-1$.
We now claim that since $G$ is a non-negative linear combination of pure
diagrams, $G'$ is as well. To see this, let $A'$ be the matrix
obtained from any diagram $A$, of codimension $p-1$, by
$$
A'_{i,j} = A_{p-1-i,p+c-j}
$$
for all integers $i$ and $j$, and note that $G'$ defined in
(\ref{GEq}) is obtained from $G$ in this way.  Remember that, by
definition, $A'$ is a diagram of codimension $p-1$, if $S_{A'}(t) =
\sum_{i,j}(-1)^i A'_{i,j}t^j$ is divisible by $(1-t)^{p-1}$. Since
$S_{A'}(t) = \sum_{i,j}(-1)^i A_{p-1-i,p+c-j}t^j = (-1)^{p-1} t^{p+c}
S_A(t^{-1})$ and this element of $\mathbb{Q}[t,t^{-1}]$ is divisible
by $(1-t)^{p-1}$ if $S_A(t)$ is, we see that $A$ is a diagram of codimension $p-1$.
Furthermore, if $\pi$ is a pure diagram of type $d =
(d_0,d_1,\dots,d_{p-1})$ then, by solving $(p-1-i,p+c-j) = (k,d_k)$
for $k=0,1,\dots,p-1$, we see that the non-zero positions of $\pi'$
are
$$
(0,p+c-d_{p-1}), (1,p+c-d_{p-1}), \dots, (p-1,p+c-d_0)
$$
and hence that $\pi'$ is a pure diagram of type
$$
d' = (p+c-d_{p-1},p+c-d_{p-2},\dots,p+c-d_0).
$$

The point of all this is that it shows that $G'$ is a non-negative
linear combination of pure diagrams and furthermore, since the maximal
shifts of $G$ are given by the type of the pure diagram
$S/\mathfrak{m}^{c-u+1}$ and this type is
$$
d = (0,c-u+1,c-u+2,\dots,c-u+p-1)
$$
we see that the minimal shifts of $G'$ are
$$
d' = (u+1,u+2,\dots,u+p-2,c).
$$
The maximal shifts of the diagram $F$ are given by the type of the
diagram $S/\mathfrak{m}^{u}$ and these are
$$
(0,u+1,u+2,\dots,u+p-1).
$$
Hence, the maximal shifts of $F$ are smaller than the minimal shifts of
$G'$, that is,
$$
\up(F) = (0,u+1,u+2,\dots,u+p-1) < (u+1,u+2,\dots,u+p-2,c) = \down(G).
$$
Since $F$ and $Q$ have the same $h$-vector, by (\ref{FCancelEq}),
$e(F) = e(Q)$ and for the same reason $e(G') = e(P^{\vee}(-c-1))$, by
(\ref{EDualEq}).  By Lemma~\ref{LefLemma}, $e(Q) = e(P^{\vee}(-c-1))$
and hence $e(F) = e(G')$. We can now apply Lemma~\ref{StepUpLemma} to
the diagrams $F$ and $G'$ which show that $E$ is a non-negative linear
combination of pure diagrams.
\end{proof}

\begin{example}
Let $R = k[x,y,z]$ and consider the $h$-vector $h(t) = 1+3t+5t^2+6t^3+2t^4$. Then, the $f_i$'s of
Proposition~\ref{LefProp} equals $f_0 = f_1 = 0$, $f_2 = 1/3$, $f_3 = 5/12$ and $f_4 = 1/4$, and
hence the diagram $F$ of the same proposition equals
\[
F =  \frac{1}{3} \beta(S/\mathfrak{m}^2)+\frac{5}{12}  \beta(S/\mathfrak{m}^3) + \frac{1}{4} \beta(S/\mathfrak{m}^4) =
\begin{pmatrix}
1 &-&- \\-& 1& 2/3\\ -& 5/3&5/4 \\ -& 5/4& 1
\end{pmatrix}.
\]
For the $g_i$'s, of Proposition~\ref{LefProp}, we get $g_0 = 0$ and $g_1 = 2$, and hence 
\[
G = 2 \cdot \beta(S/\mathfrak{m}^2)= \begin{pmatrix} 
2& -& -\\ -& 6& 4
\end{pmatrix}.
\]
The diagram $G'$, from the proof of Proposition~\ref{LefProp}, is the diagram
obtained from $G$ by rotating it 180 degrees and shifting its rows so that
its upper most non-zero element lies on row $c$, where $c$ is the degree of $h(t)$. 
In this case we get
\[
G' = \begin{pmatrix}
-&-&- \\ -&-&- \\ -&-&- \\ 4&6&-\\ -&-&2
 \end{pmatrix}.
\]
Now, the diagram $E$, of Proposition~\ref{LefProp}, defined by
$$
E_{i,j} = F_{i,j} + G'_{i-1,j} = F_{i,j} + G_{3-i,7-j},
$$ 
for all integers $i$ and $j$, is simply the sum of $F$ with the diagram obtained
from $G'$ by shifting its columns one position to the right. We get 
\[
E = \begin{pmatrix}
1 &-&-&- \\-& 1& 2/3 & -\\ -& 5/3&5/4&- \\ -& 21/4& 7&- \\ -&-&-&2
 \end{pmatrix}.
\] 
Proposition~\ref{LefProp} states that when $h(t) = 1+3t+5t^2+6t^3+2t^4$ is the
$h$-vector of an artinian level $R$-module with the weak Lefschetz property, 
then $E$ is a non-negative linear combination of
pure diagrams. For the diagram $E$ of this example this turns out to be true. We have
\begin{multline*}
E = \frac{4}{21}\pure( 0, 2, 3, 7)+\frac{1}{14}\pure ( 0, 2, 4, 7)+ \\
 \frac{13}{84}\pure( 0, 3, 4, 7)+ \frac{2}{15}\pure( 0, 3, 5, 7)+\frac{9}{20}\pure ( 0, 4, 5, 7).
\end{multline*}
In fact, Theorem~\ref{WLPUpToMultProp} says that as soon as the $f_i$'s and $g_i$'s are
non-negative, $h(t)$ is a rational multiple of the
$h$-vector of an artinian level $R$-module with the weak Lefschetz property, 
so this was no coincident.
This means that $h(t) = 1+3t+5t^2+6t^3+2t^4$ 
is a rational multiple of a level $h$-vector, 
and a non-negative linear combination of pure diagrams. Hence, by Proposition~\ref{MaxCompressedProp}, it is
a non-negative linear combination of $h$-vectors of extremely compressed diagrams. 
Computing the coefficients of these extremely compressed diagrams with Lemma~\ref{ECompLemma}
gives the diagram
\[
W = \frac{5}{21}\pure(0, 2, 3, 7) +\frac{11}{42}\pure(0, 3, 4, 7)+ \frac{1}{2}\pure(0, 4, 5, 7)=
 \begin{pmatrix}
1 &-&-&- \\-& 1& 5/6 & -\\ -& 11/6&11/6&- \\ -& 35/6& 7&- \\ -&-&-&2
 \end{pmatrix}.
\]
where $h_W(t) = 1+3t+5t^2+6t^3+2t^4$. Note that $W \geq E$, entry by entry, and hence
that $W$ could not be a rational multiple of the Betti diagram of a module with the 
weak Lefschetz property.
\end{example}

\begin{rmk}
We can replace the diagrams $F$ and $G$ of Proposition~\ref{LefProp}
with the diagrams $\MaxLex(h_Q)$ and $\MaxLex(h_{P'})$, respectively, 
where $h_Q(t)$ and $h_{P'}(t)$ are the $h$-vectors obtained from
the $h$-vector of $M$ in Lemma~\ref{LefLemma}. This gives a smaller
upper bound than the one given in  Proposition~\ref{LefProp}
but we no longer know that this upper bound
is a non-negative linear combination of pure diagrams, since we do
not know that $\MaxLex(h_Q)$ and $\MaxLex(h_{P'})$ have this property. 
In addition, the upper bound of Proposition~\ref{LefProp} is a linear function
of the $h$-vector, so when considering $h$-vectors up to multiplication
with a rational number, this upper bound behaves well, something
that is not true for one given by  $\MaxLex(h_Q)$ and $\MaxLex(h_{P'})$.
In fact, in Theorem~\ref{WLPUpToMultProp}, we will see that there always exists
an artinian $R$-module with the weak Lefschetz property having
an integer multiple of the upper bound $E$, of Proposition~\ref{LefProp},
as its Betti diagram. So when considering $h$-vectors and Betti diagrams
up to multiplication with a rational number this upper bound is in some 
sense sharp.
\end{rmk}

We will now show that given an $h$-vector such that the $f_i$'s and
$g_i$'s of Proposition~\ref{LefProp} are all non-negative, we can
construct an artinian level $R$-module with the weak Lefschetz
property whose $h$-vector is an integer multiple of this $h$-vector, and furthermore, whose Betti diagram
is an integer multiple of the maximal one given by the diagram $E$ of
Proposition~\ref{LefProp}.

To prove that our construction has the desired Betti diagram we will
use the following standard observation, and we include its simple
proof.

\begin{lemma} \label{BettiTruncLemma} Let $M$ be a graded $R$-module.
  Then, for any integer $u$,
$$
\beta_{i,j}\left(\bigoplus_{d \leq u} M_d\right) = \beta_{i,j}(M)
$$
for all $j-i < u$.
\end{lemma}

\begin{proof}
  Let
$M' = \bigoplus_{d \leq u} M_d$ and $M'' = \bigoplus_{d > u} M_d$
and consider the exact sequence
$$
0 \to M'' \to M \to M' \to 0
$$
From the associated long exact sequence of $\Tor$
$$
 \to \Tor^R_{i}(M'',k)\to \Tor^R_{i}(M,k)\to
\Tor^R_{i}(M',k) \to \Tor^R_{i+1}(M'',k)
$$
it follows that $\Tor^R_{i}(M,k)_j \cong \Tor^R_{i}(M',k)_j$ whenever
$$
\Tor^R_{i}(M'',k)_j = \Tor^R_{i+1}(M'',k)_j = 0.
$$
Since $M'_i = 0$ when $i\leq u$ we get that $\Tor^R_{i}(M',k)_j = 0$
when $j-i \leq u$ and hence that
$$
\Tor^R_{i}(M',k)_j \cong \Tor^R_{i}(M,k)_j
$$
when $j-i < u$, which gives the desired equality of the Betti numbers.
\end{proof}

By a point in $\mathbb{P}_k^{p-1}$ we mean a
$k$-rational point, and such a point $P$ can be identified with
a $p$-tuple $P = (u_0:u_1:\dots:u_{p-1})$ where $u_i \in k$.  
For any set of points $X$ denote by $I_X$ the ideal of forms vanishing on $X$, that is, the
ideal generated by 
$$
\{f \in R \,:\, f \text{ is homogeneous and } f(P)
= 0 \text{ for all } P \in X\}.
$$
The homogeneous coordinate ring of
$X$ is $R/I_X$ and we denote it by $A_X$. When $X$ is a finite set of
points, we denote the smallest integer $d$ such that $H(A_X,d) = |X|$
by $\tau(A_X)$.  

The construction in the following proposition is more general than we
need and we will only apply it to points in general position here.

\begin{prop} \label{HilbFromPointsLemma} Let $X$ be a set of distinct
  points in $\mathbb{P}^{p-1}$ and let $X = \bigcup_{i=1}^{s} Y_i$ and
  $X = \bigcup_{i=1}^{t} Z_i$ be partitions of $X$ and furthermore let
  $c$ be an integer such that $c \geq \max \{ \tau(A_{Y_i})\}_{i=1}^s
  +\max \{ \tau(A_{Z_i})\}_{i=1}^t$. Then there is an artinian
  $R$-module $M$, of codimension $p$, with the weak Lefschetz property
  and Hilbert function given by
$$
H(M,d) = \min \left\{ \sum_{i=1}^s H(A_{Y_i},d), \sum_{i=1}^t
  H(A_{Z_i},c-d) \right\}
$$
and furthermore,
$$
\beta_{i,j}(M) = \sum_{k = 1}^s \beta_{i,j}(A_{Y_k},d) + \sum_{i=1}^t
\beta_{p-1,p+c-j}(A_{Z_k},d)
$$
for all integers $i$ and $j$.
\end{prop}

\begin{proof}
  Let $M_X$ be the free $k[t,t^{-1}]$-module with basis
$\{e_P\}_{P\in X}$, that is, 
$$
M_X = \bigoplus_{P \in X} k[t,t^{-1}] e_P.
$$
We give $M_X$ the structure of a graded $R$-module by defining the
multiplication of any homogeneous element $f \in R$ and the element
$\sum_{P \in X} a_Pe_P \in M_X$, where $a_P \in k[t,t^{-1}]$, by
$$
f \cdot \sum_{P \in X} a_Pe_P = \sum_{P \in X} f(P)t^da_Pe_P
$$
where $d$ is the degree of $f$. (To make the value of $f(P)$ well
defined we have to choose a fixed coordinate representation for each
point $P \in X$, so we assume that this have been done.)
Observe that $M_X \cong \bigoplus_{d\in\mathbb{Z}}\Gamma(X,\mathcal{O}_X(d))$.

We will now construct a quotient of a submodule of $M_X$ with the
desired Hilbert function.  Let $y_{i} = \sum_{P \in Y_i} e_P$ for each
$i$, and let $M_Y$ be the submodule of $M_X$ generated by
$y_{1},y_{2},\dots,y_{s}$. Then $M_Y$ is isomorphic to the direct sum
$$
M_Y \cong \bigoplus_{i=1}^{s} Ry_{i}.
$$
We claim that $Ry_{i} \cong A_{Y_i}$, where $A_{Y_i}$ is the
homogeneous coordinate ring of $Y_i$.  Consider the $R$-module
homomorphism $R \to Ry_{i}$ given by $f \mapsto fy_{i}$, for any $f
\in R$, and note that $fy_{i}=0$ if and only if $f_0(P) = f_1(P) =
\dots = f_d(P)=0$ for all $P \in Y_i$, where $f = f_0 + f_1 + \dots +
f_d$ is the decomposition of $f$ into homogeneous polynomials. Hence,
the kernel of $R \to Ry_{i}$ is the ideal of forms vanishing on all
points of $Y_i$ which implies that $Ry_{i}$ is isomorphic to the
homogeneous coordinate ring of $Y_i$, that is, $Ry_{i} \cong A_{Y_i}$.
Since $M_Y$ is the direct sum of the $Ry_{i}$'s we get
\begin{equation} \label{HilbYEq} 
M_Y \cong \bigoplus_{i=1}^s A_{Y_i},
\end{equation}
Let $M_Z$ be the submodule generated by the elements $z_i = \sum_{P
  \in Z_i}e_P$ for $i=1,2,\dots,t$. Then we see, in the same way as we
did for $M_Y$, that
\begin{equation} \label{HilbZEq} 
M_Z \cong \bigoplus_{i=1}^t A_{Z_i}.
\end{equation}
Furthermore, let $\tau_Y = \max \{ \tau(A_{Y_i})\}_{i=1}^s$ and
$\tau_Z = \max \{ \tau(A_{Z_i})\}_{i=1}^t$. Then
\begin{equation} \label{HilbXYEq} H(M_Y,d) = |X| 
\end{equation}
for all $d\geq \tau_Y$ and
\begin{equation} \label{HilbXZEq} H(M_Z,d) =  |X|
\end{equation} 
for all $d\geq \tau_Z$ by the definition of $\tau$.

For the next step we need the $R$-module $\HHom_k(N,k) =
\bigoplus_{i\in \mathbb{Z}} \Hom_k(N_i,k)$ defined for any $R$-module
$N$. Its grading is given by $\HHom_k(N,k)_i= \Hom_k(N_{-i},k)$ for
each degree $i$, and multiplication given by letting $f\cdot \psi$,
for any $f \in R$ and $\psi \in \HHom_k(N,k)$, be the function $m
\mapsto \psi(fm)$ for any $m \in N$. We note that when $N$ is
artinian, which $M_X$, $M_Y$ and $M_Z$ are not, then $\HHom_k(N,k)
\cong N^{\vee}$, see \cite[Proposition 2.5]{Soderberg2}.

We will now show that $\HHom_k(M_X,k)(d) \cong M_X$ for any integer
$d$. For any $P \in X$ and integer $u$, let $t^{-u}e^*_P \in
\HHom_k(M_X,k)$ be defined by
$$
t^{-u}e^*_p(t^ve_Q) = \begin{cases} 1 & \text{ if $P=Q$ and $u=v$}, \\
  0 & \text{ otherwise}. \end{cases}
$$
Then for any element $f \in R$ of degree $d$ we get $f \cdot
t^{-u}e^*_P = f(P)t^{-u+d}e^*_P$ and that $(t^{-u}e^*_P)_{X \in P}$ is
a $k$-basis of $\HHom_k(M_X,k)_{-u}$, in fact, it is the $k$-dual of
the $k$-basis $(t^{u}e_P)_{P \in X}$ of $(M_X)_u$.  Hence, the
function $\HHom_k(M_X,k)(d) \to M_X$ defined by $t^{-u}e^*_P \mapsto
t^{u+d}e_P$ is an isomorphism of $R$-modules for any integer $d$.
\begin{figure}[tbp]\caption{}
  \label{HilbFig}
  \vspace{10pt} \input{hilb.pstex_t}
\end{figure}

Note that
\begin{equation} \label{MZDualEq}
H(\HHom_k(M_Z,k)(-c),d) = H(M_Z,c-d)
\end{equation}
for all $d$,
and that we seek a module, $M$, with Hilbert function given by
$H(M,d) = \max\{H(M_Y,d),H(M_Z,c-d)\}$, for all $d$. Note also
that this Hilbert function is equal to $|X|$ in degree
$\tau_Y \leq d \leq c-\tau_Z$ (see Figure~\ref{HilbFig}). 

We will now show that $\HHom_k(M_Z,k)(-c) \cong M_X/N_Z$ for some
submodule $N_Z \subseteq M_X$ and that in fact $N_Z \subseteq M_Y$.
The exact sequence
$$
0 \to M_Z \to M_X \to M_X/M_Z \to 0
$$
gives the exact sequence
$$
0 \to \HHom_k(M_X/M_Z,k) \to \HHom_k(M_X,k) \to \HHom_k(M_Z,k) \to 0.
$$
From this exact sequence and the isomorphism $\HHom_k(M_X,k)(-c) \cong
M_X$ we get
$$
0 \to \HHom_k(M_X/M_Z,k)(-c) \to M_X \to \HHom_k(M_Z,k)(-c) \to 0,
$$
and hence that $\HHom_k(M_Z,k)(-c) \cong M_X/N_Z$ for some submodule
$N_Z \subseteq M_X$. To see that $N_Z \subseteq M_Y$ observe first
that (\ref{MZDualEq}) together with the definition of $N_Z$ implies
\begin{equation} \label{HilbNZEq}
H(N_Z,d) = |X| - H(M_Z,c-d).
\end{equation}
Then it follows from the definition of $\tau_Z$ that 
the initial degree of $N_Z$ is $c-\tau_Z+1$. Finally, since
$(M_Y)_d = (M_X)_d$ for all $d \geq \tau_Y$ and $c \geq \tau_Y + \tau_Z$,
by assumption, we see that $(M_Y)_d = (M_X)_d$ for all $d \geq c-\tau_Z$, and
hence, taking into account the initial degree of $N_Z$, 
that $N_Z \subseteq M_Y$.

We now claim that the artinian $R$-module $M = M_Y/N_Z$ has the weak
Lefschetz property and the desired Hilbert function and Betti numbers.
We have the following commutative diagram with exact rows and columns
\begin{equation} \label{CDEq}
  \vcenter{\xymatrixcolsep{0.5cm}\xymatrix{
               &0 \ar[d]                    & 0 \ar[d]                    & 0 \ar[d] \\
      0 \ar[r] &N_Z \ar[r] \ar[d]^{\varphi}   &M_Y \ar[r] \ar[d]^{\varphi'}  &M \ar[d]^{\varphi''} \ar[r] &0 \\
      0 \ar[r] &\HHom(M_X/M_Z,k)(-c) \ar[r] \ar[d]  &M_X \ar[r]                    &\HHom_k(M_Z,k)(-c) \ar[r] &0 \\
               & 0 }}
\end{equation}
where $\varphi''$ is the unique homomorphism induced by
$\varphi$ and $\varphi'$.
Note that $\varphi'_d$ is an isomorphism in
degrees $d \geq \tau_Y$ and, as a consequence, $\varphi''_d$ is also
an isomorphism in degrees $d \geq \tau_Y$.

It follows from (\ref{HilbNZEq}) and the the
exact sequence in the top row of (\ref{CDEq})
that $H(M,d) = \max\{H(M_Y,d),H(M_Z,c-d)\}$ which shows that the
Hilbert function of $M$ is the desired one.
Since
$$
\bigoplus_{d \leq \tau_Y} M_d \cong \bigoplus_{d \leq \tau_Y} (M_Y)_d
$$
it follows from Lemma~\ref{BettiTruncLemma} that
\begin{equation}\label{BettiYEq}
  \beta_{i,j}(M) = \beta_{i,j}(M_Y)
\end{equation}
when $j-i < \tau_Y$.  Consider the dual of $M$, $M^{\vee}$, which
since $M$ is artinian equals $\HHom_k(M,k)$, and the homomorphism
$$
\varphi'' : M \to \HHom_k(M_Z,k)(-c)
$$
from (\ref{CDEq}) which is an isomorphism in degrees $d \geq \tau_Y$.
Applying $\HHom_k(\_,k)$ to this homomorphism and shifting the degrees
by $-c$ yields a homomorphism
$$
M_Z \to M^{\vee}(-c)
$$
which is an isomorphism in degrees $d \leq c - \tau_Y$. Hence we
have
$$
\bigoplus_{d \leq c-\tau_Y} M^{\vee}(-c)_d \cong \bigoplus_{d \leq
  c-\tau_Y} (M_Z)_d
$$
and it follows from Lemma~\ref{BettiTruncLemma} that
$$
\beta_{i,j}(M^{\vee}(-c)) = \beta_{i,j}(M_Z)
$$
when $j-i < c-\tau_Y$.  Since $\beta_{i,j}(M) =
\beta_{p-i,p+c-j}(M^{\vee}(-c))$ this means that
\begin{equation} \label{BettiZEq} \beta_{i,j}(M) =
  \beta_{p-i,p+c-j}(M_Z)
\end{equation}
when $j-i > \tau_Y$.
From (\ref{BettiYEq}) and (\ref{BettiZEq}) and the
fact that $\beta_{i,j}(M_Y) = 0$ when $j-i > \tau_Y$ and
$\beta_{p-i,p+c-j}(M_Z) = 0$ when $j-i < \tau_Y$ we get that
\begin{equation} \label{BettiUEq} \beta_{i,j}(M) = \beta_{i,j}(M_Y) +
  \beta_{p-i,p+c-j}(M_Z)
\end{equation}
for all integers $i$ and $j$ satisfying $j-i \neq \tau_Y$.  To show the
equality when $j-i = \tau_Y$, consider the diagram $D$ defined by
$$
D_{i,j} = \beta_{i,j}(M_Y) + \beta_{p-i,p+c-j}(M_Z)
$$
for all integers $i$ and $j$.  Then we have seen above that $h_D(t) =
h_M(t)$ and this implies that $S_D(t) = S_M(t)$ which gives the
equality
$$
\sum_{i,j}(-1)^i D_{i,j} t^j = \sum_{i,j}(-1)^i \beta_{i,j}(M)t^j,
$$
and if we remove all terms that are equal by (\ref{BettiUEq}) we are
left with
$$
\sum_{i=0}^p(-1)^i D_{i,u+i} t^{u+i} = \sum_{i=0}^p(-1)^i
\beta_{i,u+i}(M)t^{u+i}
$$
which shows that we have equality in (\ref{BettiUEq}) also when $j-i =
u$, and hence that $\beta(M) = D$.

It remains to show that $M$ have the weak Lefschetz property.  Let
$\ell$ be a non-zero divisor in $A_X$, that is, an element in $R$ such
that $\ell(P) \neq 0$ for all $P \in X$. Then $\ell$ is a non-zero
divisor on $M_X$ and hence also on $M_Y$ since it is a submodule of
$M_X$. This means that $(M_Y/N_Z)_{d-1} \xrightarrow{\times \ell}
(M_Y/N_Z)_{d}$ is surjective for all $d \leq c-\tau_Z$ since $(N_Z)_d
= 0$ for these values of $d$.  Assume now that $d > c-\tau_Z$.  We
have that $(M_Z)_{c-d-1} \xrightarrow{\times \ell} (M_Z)_{c-d}$ is
surjective, since $M_Z$ is a submodule of $M_X$, and by applying
$\HHom_k(\_,k)$ we get that $\HHom_k((M_Z)_{c-d},k) \xrightarrow{\ell}
\HHom_k((M_Z)_{c-d-1},k)$ is injective.  Since $(M_Y/N_Z)_d \cong
(M_X/N_Z)_d \cong \HHom_k((M_Z)_{c-d},k)$, when $d > c-\tau_Z$, we get
that $(M_Y/N_Z)_d \xrightarrow{\ell} (M_Y/N_Z)_{d+1}$ is injective
which finishes the proof.
\end{proof}

\begin{thm} \label{WLPUpToMultProp} Assume that $k$ is an infinite
  field and that $h(t) = \sum_{i\in \mathbb{Z}}h_it^i$ is a polynomial
  of degree $c$ and furthermore, let $u$ be the smallest integer such
  that $h_u \geq h_{u+1}$. Then $h(t)$ is a rational multiple of the
  $h$-vector of a artinian level $R$-module of codimension $p$,
  generated in degree zero, having the weak Lefschetz property, if and only if
$$
f_i = \frac{h_i-h_{i-1}}{s_i} - \frac{h_{i+1} - h_{i}}{s_{i+1}} \geq 0,
$$
for $i= 0,1,\dots,u-1$, and
$$
g_i = \frac{h_{c-i} - h_{c-i+1}}{s_i} - \frac{h_{c-i-1} -
  h_{c-i}}{s_{i+1}} \geq 0,
$$
for $i = 0,1,\dots,c-u-1$, where $s_i = \binom{p-2+i}{p-2}$. Furthermore, if the $f_i$'s and
$g_i$'s are all non-negative, then there exist an artinian level
$R$-module with the weak Lefschetz property whose Betti diagram
is an integer multiple of the diagram $E$ of Proposition~\ref{LefProp}.
\end{thm}

\begin{proof}
  The necessity of the conditions follows immediately from
  Proposition~\ref{LefProp}.

  Assume that the $f_i$'s and $g_i$'s of the proposition are all
  non-negative and furthermore, let
$$
f_u = \frac{h_u-h_{u-1}}{s_u}
$$
and
$$
g_{c-u} =\frac{h_{u} - h_{u+1}}{s_{c-u}}.
$$
Then $f_u \geq 0$ and $g_{c-u} \geq 0$ since, by assumption, $u$ is
the smallest integer such that $h_u \geq h_{u+1}$.  Let $m$ be an
integer such that $mf_i$ is an integer for $i=0,1,\dots,u$ and $mg_i$
is an integer for $i=0,1,\dots,c+u$.  We will use
Proposition~\ref{HilbFromPointsLemma} to show the existence
of an artinian level $R$-module with the weak Lefschetz property
and Betti diagram the diagram $mE$, for some integer $m$, where $E$ is the
diagram of Proposition~\ref{LefProp}.
Then $$
h(t) = h_E(t) = \frac{1}{m}h_M(t)
$$ 
which shows that $h(t)$ is a rational multiple of the $h$-vector of $M$.

Consider a set, $X$, of points in $\mathbb{P}^{p-1}$ in general
position.  It is a fact that, if $k$ is an infinite field, the
condition on $X$ that makes the Hilbert function of $A_X$ given by
$H(A_X,d) = \min \{r_d,|X|\}$ is open and non-empty.
 As a consequence,
if we also require that $H(A_Y,d) = \min\{r_d,|Y|\}$ for any subset
$Y$ of $X$ the condition is still non-empty and open condition. Assume
that $X$ satisfies this condition and that there are two partitions
$$
X = \bigcup_{j=0}^u \bigcup_{i=1}^{mf_j} Y_{i,j} \quad\text{ and
}\quad X = \bigcup_{j=0}^{c-u} \bigcup_{i=1}^{mg_j} Z_{i,j}
$$
such that $|Y_{i,j}| = |Z_{i,j}| = r_j$ for all $i$ and $j$, where
$r_j = \binom{p-1+j}{p-1}$.  Then $H(A_{Y_{i,j}},d) =
\min\{r_d,|Y_{i,j}|\} = \min\{r_d, r_j\}$ and $H(A_{Z_{i,j}},d) =
\min\{r_d,|Z_{i,j}|\} = \min\{r_d,r_j\}$ for each integer $d$. 

Since $\tau(A_{Y_{i,j}}) = \tau(A_{Z_{i,j}}) = j$,
for each integer $j$, we get 
\begin{multline*}
\max\{\tau(A_{Y_{i,j}}) \,|\,0\leq j \leq u, 1 \leq i \leq mf_j \} \\
+ \max\{\tau(A_{Z_{i,j}}) \,|\,0\leq j \leq c-u, 1 \leq i \leq mg_j \} \leq u + c - u = c.
\end{multline*}
Note that the above is an inequality and not an equality since
 $g_j$ might be zero
for some integers $j$.
Proposition~\ref{HilbFromPointsLemma} now shows
the existence of an
artinian $R$-module, $M$, with the weak Lefschetz property and 
Betti diagram given by
\begin{equation} \label{BettiLef}
\beta_{k,l}(M) = 
\sum_{j=0}^u\sum_{i=1}^{mf_j} \beta_{k,l}(A_{Y_{i,j}})
+\sum_{j=0}^{c-u}\sum_{i=1}^{mg_j} \beta_{p-k,p+c-l}(A_{Z_{i,j}})
\end{equation}
By the way the $A_{Y_{i,j}}$'s and
$A_{Z_{i,j}}$'s are defined we see that
$$
\beta^R(A_{Y_{i,j}}) = \beta^S(S/\mathfrak{m}^{j+1})
$$
and
$$
\beta^R(A_{Z_{i,j}}) = \beta^S(S/\mathfrak{m}^{j+1})
$$
Hence
$$
\sum_{j=0}^u\sum_{i=1}^{mf_j} \beta(A_{Y_{i,j}})  = m\sum_{i=0}^uf_i\beta^S(S/\mathfrak{m}^{j+1})
$$
and we see that this equals $mF$, where $F$ is as defined in
Proposition~\ref{LefProp}. In the same way we get
$$
\sum_{j=0}^{c-u}\sum_{i=1}^{mg_j} \beta(A_{Z_{i,j}})  = m\sum_{i=0}^{c-u}g_i\beta^S(S/\mathfrak{m}^{j+1})
$$
which then equals $mG$, where $G$ is as defined in
Proposition~\ref{LefProp}. And since, by (\ref{BettiLef}), $\beta(M)$ then satisfies
$$
\beta_{i,j}(M) = mF_{i,j} + mG_{p-i,p+c-j} = mE_{i,j}.
$$
where $E$ is as defined in Proposition~\ref{LefProp}. 
We have shown that $h(t)$ is a rational multiple of the
$h$-vector of a level artinian $R$-module with the weak Lefschetz
property whose Betti diagram is the diagram $E$ of
Proposition~\ref{LefProp}.
\end{proof}

\begin{rmk}
  In codimension two the conditions on the $f_i$'s and $g_i$'s of
  Theorem~\ref{WLPUpToMultProp} are equivalent to
$$
h_{i+1} - 2h_i + h_{i-1} \leq 0.
$$
for $i=0,1,\dots,c$.  This follows by straightforward calculation
using $s_i = 1$.  This in turn is equivalent to the condition of
Theorem~\ref{NormHVec}, in the codimension two case, and
furthermore, this precisely describes the level $h$-vectors of modules
of codimension two (see Remark~\ref{LevelCodim2Rmk}).
\end{rmk}

\section{Level modules of codimension three} \label{CodimThree}
We prove a result on linear combinations of pure diagrams of 
codimension three. Combining this result with the result
on Lefschetz modules we can prove that the Betti diagram of
any level $R$-module
of codimension three whose artinian reduction has the weak
Lefschetz property is a non-negative linear combination
of pure diagrams. This in turn proves the Multiplicity
conjecture for these module and even a stronger
conjecture of Zanello. We finish this section with a proof
of the upper bound of the Multiplicity conjecture of
Herzog, Huneke and Srinivasan for
level modules of codimension three, with or without the weak Lefschetz property. 

Zanello has proposed a strengthening of the
Multiplicity conjecture in the case of level algebras of codimension
three.  The conjecture of Zanello \cite[Conj.
2.3]{Migliore-Nagel-Zanello} is that the Multiplicity conjecture is
true not only for the Betti diagram of any level algebra of
codimension three but also for any diagram obtained from that diagram
by doing all possible cancellations (for level algebras of codimension
three this maximally cancelled, and hence in a sense minimal, diagram
is unique).  This can be formulated in terms of the $h$-vector of the
level algebra $R/I$, since the minimal and maximal shifts for the
maximally cancelled diagram can be obtained from the polynomial
$S_{R/I}(t) = \sum_{i,j} (-1)^i \beta_{i,j}(R/I)t^j$ and, since $R/I$
is Cohen-Macaulay of codimension three, $S_{R/I}(t) = (1-t)^3h(t)$.
Assume that $h(t)$ is of degree $c$ and 
$S_{R/I}(t) = \sum_{i=0}^{c+3}\Delta^3h_it^i$ and let
\begin{enumerate}
\item[\emph{i)}]  $n_1 = \min\{i \;|\; \Delta^3h_i < 0\}$, 
\item[\emph{ii)}] $n_2 = \min\{i > 0\;|\; \Delta^3h_i > 0 \}$,
\item[\emph{iii)}] $N_1 = \max\{i \leq c+1\;|\; \Delta^3h_i < 0\}$ and
\item[\emph{iv)}] $N_2 = \max\{i\;|\;  \Delta^3h_i > 0\}$,
\end{enumerate}
that is, $n_1$ and 
$n_2$ are the degrees of the first and second sign
changes in $S_{R/I}(t)$, respectively, and $N_1$ and $N_2$ the degrees
of the second last and last sign changes.

\begin{conj}[Zanello] \label{Zanello} If $R/I$ is level of codimension
  three then
$$
\frac{1}{3!} n_1n_2(c+3) \leq e(R/I) \leq \frac{1}{3!} N_1N_2(c+3).
$$
\end{conj}
Conjecture~\ref{Zanello} has been proved for Gorenstein algebras of
codimension three by Migliore, Nagel and Zanello
\cite{Migliore-Nagel-Zanello}.

We will now show that Conjecture~\ref{Zanello} holds for any level
algebra of codimension three whose Betti diagram is a non-negative
linear combination of pure diagrams.  The reason for this is that this
set of diagrams is, as we will see, closed under the operation
of making cancellations. Hence the strengthening proposed by Zanello
follows from Conjecture~\ref{MainConj}.

We also note that since the Betti diagram of a Gorenstein codimension
three algebra by \cite[Theorem 4.3]{Boij-Soderberg} is a non-negative
linear combination of pure diagrams we get an alternative proof, using
the technique of pure diagrams, of Zanello's stronger conjecture in
this case.

\begin{prop} \label{Codim3CancelProp} Let $D$ be a non-negative linear
  combination of pure diagrams of codimension three. If $D$ has only
  one non-zero entry in column zero and column three, then any
  consecutive cancellation of $D$
 is still a non-negative linear combination of pure diagrams
\end{prop}

\begin{proof}
  Note that, in codimension one, any consecutive cancellation
of a non-negative linear combination
of pure diagrams is still a non-negative linear combination
of pure diagrams.
In fact, let $D$ be a diagram of codimension one.
  Since for any shifts $(d_0,d_1)$ the pure diagram $\pure(d_0,d_1)$
  has as its non-zero entries $\pure(d_0,d_1)_{0,d_0}=1$ and
  $\pure(d_0,d_1)_{1,d_1}=1$ we see that $D$ is a convex combination
  of pure diagrams if and only if all entries of $D$ are non-negative
  and
  \begin{equation*} \label{eqPDimOne_2} \sum_{i \leq l} D_{0,i} \geq
    \sum_{i \leq l} D_{1,i+1}
  \end{equation*}
  for each integer $l$, and furthermore, we see that any cancellation
  of $D$ still satisfies these conditions.

  Now, denote by $\up = (\up_0,\up_1,\up_2,\up_3)$ and $\down =
  (\down_0,\down_1,\down_2,\down_3)$ the maximal and minimal shifts of
  $D$, respectively.  Then, by the hypothesis on the entries of $D$,
  we have $\up_0 = \down_0$ and $\up_3 = \down_3$. This permits us to
  use the isomorphisms $\phi_0$ and $\phi_3$ from
  Definition~\ref{stepDown}. We get an isomorphism
$$
\phi_0 \circ \phi_3 :\, \vecsp \to
V_{(\down_1,\down_2),(\up_1,\up_2)}.
$$
Since this isomorphism preserves non-negative linear combination of pure diagrams
the proposition follows from the codimension one case above.
\end{proof}

\begin{cor} \label{Codim3CancelCor} 
  Let $M$ be a level module of codimension three whose Betti diagram, $\beta(M)$,
is a non-negative linear combination of pure diagrams, then
any diagram obtained from $\beta(M)$ 
by a sequence of consecutive cancellations is a non-negative linear combination of pure diagrams. 
\end{cor}

\begin{proof}
Since $M$ is level, $\beta(M)$ has only one non-zero entry in column
zero and column three, and hence we can apply Proposition~\ref{Codim3CancelProp}.
\end{proof}

\begin{cor} \label{ZanelloCodim3Cor}  Conjecture~\ref{Zanello}, and
hence the Multiplicity conjecture, is true
  for any codimension three level algebra whose Betti diagram is a
  non-negative linear combination of pure diagrams. 
\end{cor}

\begin{proof}
Let $M$ be a level algebra with a Betti diagram, $\beta(M)$, which is a non-negative
linear combination of pure diagrams.
The unique diagram, $D$, obtained
from $\beta(M)$ by doing all possible consecutive cancellations,
is then still a non-negative linear combination of pure diagrams, by Corollary~\ref{Codim3CancelCor}.
Then the condition of the Multiplicity conjecture hold for $D$,
since by \cite[Proposition 2.8]{Boij-Soderberg} it
holds for any non-negative linear combination of pure diagrams, and this condition
for the diagram $D$ is precisely the condition of Conjecture~\ref{Zanello} for $\beta(M)$.
\end{proof}

\begin{prop} \label{CodimThreeWLPCor} Let $M$ be a level $R$-module
  whose artinian reduction has the weak Lefschetz property. If $M$
  is of codimension three, then its Betti diagram is a
  non-negative linear combination of pure diagrams which implies
$$
\beta_0(M)\frac{\down_1\down_2\down_3}{3!} \leq e(M) \leq
\beta_0(M)\frac{\up_1\up_2\up_3}{3!}.
$$
where $(0,\down_1,\down_2,\down_3)$ and $(0,\up_1,\up_2,\up_3)$ are
the minimal and maximal shifts of $M$, respectively. Moreover, $M$
even satisfies the stronger condition of Zanello's
Conjecture~\ref{Zanello}.
\end{prop}

\begin{proof}
By Proposition~\ref{LefProp}, $\beta(M)$ is obtained
from a non-negative linear combination of pure diagrams
by a sequence of consecutive cancellations. Hence $\beta(M)$
is non-negative linear combination of pure diagrams
by Corollary~\ref{Codim3CancelCor}. The two other assertions
follows from Corollary~\ref{ZanelloCodim3Cor}.
\end{proof}

\subsection{The upper bound of the Multiplicity conjecture}
In this section we show that level modules of codimension three
satisfies the upper bound of the Multiplicity conjecture of Herzog,
Huneke and Srinivasan. The technique we use is based on the
description, given in Proposition~\ref{NormMaxBettiProp}, of Betti diagrams as
cancellations of non-negative linear combinations of pure diagrams.

We need the following result from \cite[Corollary
2.12]{Boij-Soderberg}.

\begin{prop} \label{QuasiPureProp} If $E$ is a diagram of codimension
  $\pdim$ with maximal and minimal shifts $\up =
  (\up_0,\up_1,\dots,\up_\pdim)$ and $\down =
  (\down_0,\down_1,\dots,\down_\pdim)$ satisfying $\up_{i-1} <
  \down_i$ for all $0 \leq i \leq \pdim$ then $E$ is a non-negative
  linear combination of pure diagrams.
\end{prop}

Note that proposition~\ref{QuasiPureProp} is similar to a result by Herzog
and Srinivasan, generalized to modules by Migliore, Nagel and
R\"{o}mer \cite{Migliore-Nagel-Roemer}, stating that any module with a
quasipure minimal free resolution, that is, with maximal and minimal
shifts satisfying $\up_{i-1} \leq \down_i$, satisfies the Multiplicity
conjecture. 

The proof of Theorem~\ref{MCUpperThm} uses the following idea from
\cite{Boij-Soderberg}.  The multiplicity of the pure diagram
$\pure(0,d_1,d_2,\dots,d_p)$ is $e(\pure(d)) = d_1d_2 \cdots d_p/p!$.
Hence, for any diagram $D = \sum_{d^{\prime} \leq d}
a_{d^{\prime}}\pure(d^{\prime})$, where $a_{d^{\prime}}$ are
non-negative integers, we get
\begin{equation} \label{MCFromCombEq}
e(D) = \sum_{d^{\prime} \leq d} a_{d^{\prime}}e(\pure(d^{\prime}))
\leq \sum_{d^{\prime} \leq d} a_{d^{\prime}}e(\pure(d)) = D_{0,0}
\frac{d_1d_2\cdots d_p}{p!}.
\end{equation}

\begin{thm} \label{MCUpperThm} Any level $R$-module, $M$, of
  codimension three, with maximal shifts given by $\up =
  (0,\up_1,\up_2,\up_3)$, satisfies
$$
e(M) \leq \beta_0(M) \frac{\up_1 \up_2 \up_3}{3!}.
$$
\end{thm}

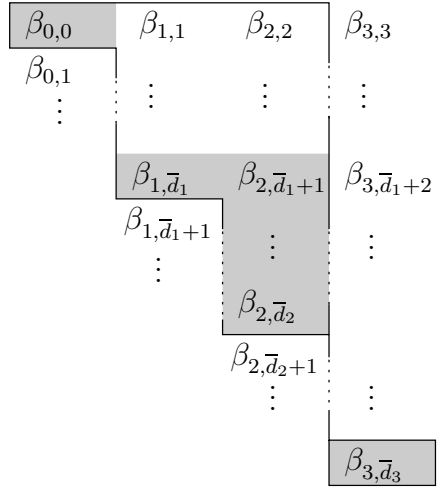
\begin{figure}[tbp]\caption{The Betti diagram of a level module $M$
    with maximal shifts $\up = (0,\up_1,\up_2,\up_3)$. Every Betti
    number outside the outlined area are zero. The grey area indicates the
    non-zero entries of the diagram $E$ in the proof of
    Theorem~\ref{MCUpperThm}.}
  \label{BettiFig}
  \vspace{10pt} \input{bettidiagram.pstex_t}
\end{figure}

\begin{proof}
  Let $M$ be a level $R$-module of codimension three, with maximal
  shifts given by $\up = (0,\up_1,\up_2,\up_3)$. By artinian
reduction we may assume that $M$ is artinian and hence that
$R = k[x_1,x_2,x_3]$.
Note that $\beta(M)$
in degrees between $\up_1+1$ and $\up_2$ have non-zero entries
only in column two, see Figure~\ref{BettiFig}. We will now show
that there is a diagram $F=E+D'$ such that $h_F(t) = h_M(t)$
and where $E$ have non-zero entries only in the positions
indicated by the grey area in Figure~\ref{BettiFig}, and
$D'$ have maximal shifts $(0,d_1-2,d_1-1,d_1)$. 

 We start with the
  description of $\beta(M)$ given in
  Proposition~\ref{NormMaxBettiProp}. Since 
$\beta(R/\mathfrak{m}^{j+1})$ is the pure diagram $\pure(0,j+1,j+2,j+3)$
we have, by Proposition~\ref{NormMaxBettiProp},
$$
\beta(M) = D - \sum_{i,j}
b_{i,j}C^{i,j}
$$
where
$$
D = \sum_{j=0}^c a_j \pure(0,j+1,j+2,j+3)
$$
for some non-negative rational numbers $a_j$ and
$b_{i,j}$ and $c = \up_3-3$.  Consider the diagrams 
$$
D' = \sum_{j=0}^{\up_1-2} a_j\pure(0,j+1,j+2,j+3)
$$
and
$$
D'' = \sum_{j=\up_1-1}^c a_j\pure(0,j+1,j+2,j+3)
$$
and note that $D = D' + D''$.

We will now see that the diagram, $E$, mentioned in the beginning of this,
can be obtained from $D''$ by a sequence of consecutive cancellations.
Note that $D''$ is zero in all positions
$(i,j)$ such that $j-i < \up_1-1$ and that the only non-zero entries in the
matrix $C^{i,j}$ are $C^{i,j}_{i,j} = 1$ and $C^{i,j}_{i+1,j} = 1$.
Now we want to take as many terms $\beta_{i,j}C^{i,j}$ from
$\sum_{i,j} b_{i,j}C^{i,j}$ and remove them from $D''$ without getting
any negative entries in the result.  If we do this we end up with
$$
E = D'' - \sum_{j=\up_1}^c \left( b_{1,j+1}C^{1,j+1} +
  b_{2,j+2}C^{2,j+2} \right).
$$
Note that $\beta_{i,j}(M) = E_{i,j}$ whenever $j-i \geq \up_1-1$ or whenever $i = 3$.
The maximal shifts of $E$ are thus the same as those of $\beta(M)$ and
$E$ has non-zero entries only in the positions indicated by the grey area
in Figure~\ref{BettiFig}.
By Proposition~\ref{QuasiPureProp}  we get that $E$ is a non-negative
linear combination of pure diagrams $\pure(d)$ whose type satisfies $d
\leq (0,\up_1,\up_2,\up_3)$ and hence the same is true for the diagram
$$
F = E + D'.
$$
The diagram $F$ have the same $h$-vector as $M$, since the cancellations does not
affect the $h$-vector, and since $F$ is a non-negative linear combination
of pure diagrams $\pure(d)$ whose type satisfies $d \leq
(0,\up_1,\up_2,\up_3)$ we get, by (\ref{MCFromCombEq}),
$$
e(M) = h_M(1) = h_F(1) = e(F) \leq \beta_0(M) \frac{\up_1 \up_2
  \up_3}{3!}.
$$
\end{proof}

What makes the proof of Theorem~\ref{MCUpperThm} possible is that the
upper bound of the Multiplicity conjecture is true for any
cancellation of a diagram on the form
$$
\sum_{j=0}^c a_j \pure(0,j+1,j+2,j+3)
$$
for some non-negative integers $a_0,a_1,\dots,a_c$, as long as it has
only one entry in column three. This means that we can prove
Theorem~\ref{MCUpperThm} without knowing which cancellations of this
diagram that actually are Betti diagrams of some module.  To see that
this is not always the case, we now give two examples.  The first,
Example \ref{LowerBoundExample}, shows that the lower bound of
Multiplicity conjecture can not be shown to hold with this technique.
The second, Example \ref{NonLevelExample}, shows that we need the
module to be level (actually it is enough that all but one of the
entries in the last column of its Betti diagram can be numerically
cancelled).
 
\begin{example} \label{LowerBoundExample} Consider the $h$-vector
  $h(t) = 16 + 48t + 21t^2 + 10t^3$. The diagram $D = \sum_{j=0}^3 a_j
  \pure(0,j+1,j+2,j+3)$, where $a_j = h_j/r_j - h_{j+1}/r_{j+1}$, is then
$$
D = \begin{pmatrix} 16 &-&-&-\\ -&75 & 100& 75/2\\
  -&25&75/2&15\\-&15&24&10 \end{pmatrix}.
$$
This diagram can be cancelled, by consecutive cancellations, to the diagram
$$
\begin{pmatrix} 16 &-&-&-\\ -&75 & 75&- \\ -&-&-&-\\-&15&9&10
\end{pmatrix}
$$
with minimal shifts $\down_1 = 2$, $\down_2 = 3$ and $\down_3 = 6$.
The lower bound of the Multiplicity conjecture for these shifts are
$$
16 \frac{\down_1 \down_2 \down_3}{6} =16\frac{2 \cdot 3 \cdot 6}{6} =
96,
$$
while its multiplicity is $h(1) = 16+48+21+10=95$.
\end{example}

\begin{example} \label{NonLevelExample} Now consider the $h$-vector
  $h(t) = 5+15t+18t^2+15t^3$. The diagram $D = \sum_{j=0}^3 a_j
  \pure(0,j+1,j+2,j+3)$, where $a_j = h_j/r_j - h_{j+1}/r_{j+1}$, is then
$$
D = \begin{pmatrix} 5&-&-&-\\-&12&16&6\\-&15&45/2&9\\-&45/2&36&15
\end{pmatrix}.
$$
This diagram can be cancelled, by consecutive cancellations, to the diagram
$$
\begin{pmatrix} 5&-&-&-\\-&12&1&6\\-&-&-&-\\-&-&27&15 \end{pmatrix}.
$$
with maximal shifts $\up_1 = 2$, $\up_2 = 5$ and $\up_3 = 6$. The
upper bound of the Multiplicity conjecture for these shifts are
$$
5 \frac{\up_1 \up_2 \up_3}{6} = 5\frac{2 \cdot 5 \cdot 6}{6}= 50,
$$
while its multiplicity is $5+15+18+15=53$.
\end{example}

\begin{rmk}
To prove the lower bound of the Multiplicity conjecture it is enough to consider
level modules. In fact, assume that $M$ is an artinian $R$-module of codimension $p$, 
generated in degree zero,
with minimal shifts given by $\down = (0,\down_1,\down_2,\dots,\down_p)$.
Let $M'$ be the $R$-module
$$
M' = \bigoplus_{i = 0}^{d_p-p} M_i
$$
The Betti numbers of this module satisfies
$$ 
\beta_{i,j}(M') = \beta_{i,j}(M),
$$
for all $j-i < d_p-p$, by Lemma~\ref{BettiTruncLemma}. This means that
the minimal shifts of $M'$ and $M$ are the same and since furthermore
$$
e(M) = \sum_i \dim_k M_i \geq e(M') = \sum_{i\leq d_p-p} M_i
$$
we see that if the lower bound of the Multiplicity conjecture holds for
$M'$ it holds for $M$ as well, that is,
$$
\beta_0(M) \frac{\down_1 \down_2 \dots \down_p}{p!} \leq e(M') \leq e(M).
$$
\end{rmk}

\subsection{Convexity of the order complex}
In \cite{Boij-Soderberg} we introduced the order complex
$\Delta(\poset)$ associated to a pair of strictly increasing sequences
$\down = (\down_0,\down_1,\dots,\down_p)$ and $\up =
(\up_0,\up_1,\dots,\down_p)$.  Consider the partial order on strictly
increasing sequences of length $p+1$ defined by $d \leq d^{\prime}$ if
$d_i \leq d^{\prime}_i$ for each $i = 0,1,\dots,p$.  The order complex
$\Delta(\poset)$ is a geometric realization of the the order complex
associated to the partial ordered set of all strictly increasing
sequences $d$ such that $\down \leq d \leq \up$.

Boij and the author conjectured $\Delta(\poset)$ to be a convex set,
and noted that this is equivalent to the fact that any non-negative
linear combination of pure diagrams can be written as a non-negative
linear combination of pure diagrams from the same chain (where the
pure diagrams inherits the order from their types). We will now show
that this is true for level algebras of codimension three.

\begin{prop}
  Let $D$ be a non-negative linear combination of pure diagrams of
  codimension three. If $D$ has only one non-zero entry in
  column zero and column three, then $D$ is a non-negative linear
  combination of pure diagrams all from the same chain.  In other
  words, if $\down = (\down_0,\down_1,\down_2,\down_3)$ and $\up =
  (\up_0,\up_1,\up_2,\down_3)$ are strictly increasing sequences of
  integers such that $\down_0 = \up_0$ and $\down_3 = \up_3$, then
  $\Delta(\poset)$ is a convex set.
\end{prop}

\begin{proof}
  By \cite[Theorem 3.4]{Boij-Soderberg}, $\Delta(\varPi_{\down',\up'})$ where 
$\down' = (\down_0,\down_1,\down_2)$ and $\up' = (\up_0,\up_1,\up_2)$
is a convex set. The isomorphism
$$
\frac{1}{\up_3}\phi_3 : V_{\down,\up} \to V_{\down',\up'}
$$
of Definition~\ref{stepDown} preserves convex combinations, and furthermore
the image of $\Delta(\poset)$ under $\frac{1}{\up_3}\phi_3$ is $\Delta(\varPi_{\down',\up'})$.
This shows that $\Delta(\poset)$ is convex.
\end{proof}

\bibliographystyle{abbrv} \bibliography{algbib}
\end{document}

%% file: hilb.pstex_t
\begin{picture}(0,0)%
\includegraphics{hilb.pstex}%
\end{picture}%
\setlength{\unitlength}{4144sp}%
\begingroup\makeatletter\ifx\SetFigFont\undefined%
\gdef\SetFigFont#1#2#3#4#5{%
  \reset@font\fontsize{#1}{#2pt}%
  \fontfamily{#3}\fontseries{#4}\fontshape{#5}%
  \selectfont}%
\fi\endgroup%
\begin{picture}(3907,2464)(-11,601)
\put( 91,659){\makebox(0,0)[lb]{\smash{{\SetFigFont{12}{14.4}{\rmdefault}{\mddefault}{\updefault}{\color[rgb]{0,0,0}$0$}%
}}}}
\put(991,659){\makebox(0,0)[lb]{\smash{{\SetFigFont{12}{14.4}{\rmdefault}{\mddefault}{\updefault}{\color[rgb]{0,0,0}$\tau_Y$}%
}}}}
\put(1531,659){\makebox(0,0)[lb]{\smash{{\SetFigFont{12}{14.4}{\rmdefault}{\mddefault}{\updefault}{\color[rgb]{0,0,0}$c-\tau_Z$}%
}}}}
\put(2881,659){\makebox(0,0)[lb]{\smash{{\SetFigFont{12}{14.4}{\rmdefault}{\mddefault}{\updefault}{\color[rgb]{0,0,0}$c$}%
}}}}
\put(2701,1109){\makebox(0,0)[lb]{\smash{{\SetFigFont{12}{14.4}{\rmdefault}{\mddefault}{\updefault}{\color[rgb]{0,0,0}$H(M,d)$}%
}}}}
\put(2611,1919){\makebox(0,0)[lb]{\smash{{\SetFigFont{12}{14.4}{\rmdefault}{\mddefault}{\updefault}{\color[rgb]{0,0,0}$H(M_Z,c-d)$}%
}}}}
\put(2881,2909){\makebox(0,0)[lb]{\smash{{\SetFigFont{12}{14.4}{\rmdefault}{\mddefault}{\updefault}{\color[rgb]{0,0,0}$H(M_Y,d)$}%
}}}}
\put( 91,1469){\makebox(0,0)[lb]{\smash{{\SetFigFont{12}{14.4}{\rmdefault}{\mddefault}{\updefault}{\color[rgb]{0,0,0}$0$}%
}}}}
\put( 91,2279){\makebox(0,0)[lb]{\smash{{\SetFigFont{12}{14.4}{\rmdefault}{\mddefault}{\updefault}{\color[rgb]{0,0,0}$0$}%
}}}}
\put(991,2279){\makebox(0,0)[lb]{\smash{{\SetFigFont{12}{14.4}{\rmdefault}{\mddefault}{\updefault}{\color[rgb]{0,0,0}$\tau_Y$}%
}}}}
\put(1531,1469){\makebox(0,0)[lb]{\smash{{\SetFigFont{12}{14.4}{\rmdefault}{\mddefault}{\updefault}{\color[rgb]{0,0,0}$c-\tau_Z$}%
}}}}
\put(2881,1469){\makebox(0,0)[lb]{\smash{{\SetFigFont{12}{14.4}{\rmdefault}{\mddefault}{\updefault}{\color[rgb]{0,0,0}$c$}%
}}}}
\put(3376,794){\makebox(0,0)[lb]{\smash{{\SetFigFont{12}{14.4}{\rmdefault}{\mddefault}{\updefault}{\color[rgb]{0,0,0}$d$}%
}}}}
\put(3376,1604){\makebox(0,0)[lb]{\smash{{\SetFigFont{12}{14.4}{\rmdefault}{\mddefault}{\updefault}{\color[rgb]{0,0,0}$d$}%
}}}}
\put(3376,2414){\makebox(0,0)[lb]{\smash{{\SetFigFont{12}{14.4}{\rmdefault}{\mddefault}{\updefault}{\color[rgb]{0,0,0}$d$}%
}}}}
\end{picture}%

%% file: bettidiagram.pstex_t
\begin{picture}(0,0)%
\includegraphics{bettidiagram.pstex}%
\end{picture}%
\setlength{\unitlength}{4144sp}%
\begingroup\makeatletter\ifx\SetFigFont\undefined%
\gdef\SetFigFont#1#2#3#4#5{%
  \reset@font\fontsize{#1}{#2pt}%
  \fontfamily{#3}\fontseries{#4}\fontshape{#5}%
  \selectfont}%
\fi\endgroup%
\begin{picture}(2544,2904)(4039,-4123)
\put(6031,-4021){\makebox(0,0)[lb]{\smash{{\SetFigFont{12}{14.4}{\rmdefault}{\mddefault}{\updefault}{\color[rgb]{0,0,0}$\beta_{3,\up_3}$}%
}}}}
\put(5356,-3391){\makebox(0,0)[lb]{\smash{{\SetFigFont{12}{14.4}{\rmdefault}{\mddefault}{\updefault}{\color[rgb]{0,0,0}$\beta_{2,\up_2+1}$}%
}}}}
\put(4726,-2581){\makebox(0,0)[lb]{\smash{{\SetFigFont{12}{14.4}{\rmdefault}{\mddefault}{\updefault}{\color[rgb]{0,0,0}$\beta_{1,\up_1+1}$}%
}}}}
\put(5401,-3076){\makebox(0,0)[lb]{\smash{{\SetFigFont{12}{14.4}{\rmdefault}{\mddefault}{\updefault}{\color[rgb]{0,0,0}$\beta_{2,\up_2}$}%
}}}}
\put(5446,-1411){\makebox(0,0)[lb]{\smash{{\SetFigFont{12}{14.4}{\rmdefault}{\mddefault}{\updefault}{\color[rgb]{0,0,0}$\beta_{2,2}$}%
}}}}
\put(4141,-1411){\makebox(0,0)[lb]{\smash{{\SetFigFont{12}{14.4}{\rmdefault}{\mddefault}{\updefault}{\color[rgb]{0,0,0}$\beta_{0,0}$}%
}}}}
\put(4816,-1411){\makebox(0,0)[lb]{\smash{{\SetFigFont{12}{14.4}{\rmdefault}{\mddefault}{\updefault}{\color[rgb]{0,0,0}$\beta_{1,1}$}%
}}}}
\put(4906,-2896){\makebox(0,0)[lb]{\smash{{\SetFigFont{12}{14.4}{\rmdefault}{\mddefault}{\updefault}{\color[rgb]{0,0,0}$\vdots$}%
}}}}
\put(5581,-3661){\makebox(0,0)[lb]{\smash{{\SetFigFont{12}{14.4}{\rmdefault}{\mddefault}{\updefault}{\color[rgb]{0,0,0}$\vdots$}%
}}}}
\put(6166,-3661){\makebox(0,0)[lb]{\smash{{\SetFigFont{12}{14.4}{\rmdefault}{\mddefault}{\updefault}{\color[rgb]{0,0,0}$\vdots$}%
}}}}
\put(5581,-2761){\makebox(0,0)[lb]{\smash{{\SetFigFont{12}{14.4}{\rmdefault}{\mddefault}{\updefault}{\color[rgb]{0,0,0}$\vdots$}%
}}}}
\put(6166,-2761){\makebox(0,0)[lb]{\smash{{\SetFigFont{12}{14.4}{\rmdefault}{\mddefault}{\updefault}{\color[rgb]{0,0,0}$\vdots$}%
}}}}
\put(4861,-1861){\makebox(0,0)[lb]{\smash{{\SetFigFont{12}{14.4}{\rmdefault}{\mddefault}{\updefault}{\color[rgb]{0,0,0}$\vdots$}%
}}}}
\put(5536,-1861){\makebox(0,0)[lb]{\smash{{\SetFigFont{12}{14.4}{\rmdefault}{\mddefault}{\updefault}{\color[rgb]{0,0,0}$\vdots$}%
}}}}
\put(6121,-1861){\makebox(0,0)[lb]{\smash{{\SetFigFont{12}{14.4}{\rmdefault}{\mddefault}{\updefault}{\color[rgb]{0,0,0}$\vdots$}%
}}}}
\put(6031,-1411){\makebox(0,0)[lb]{\smash{{\SetFigFont{12}{14.4}{\rmdefault}{\mddefault}{\updefault}{\color[rgb]{0,0,0}$\beta_{3,3}$}%
}}}}
\put(4321,-1951){\makebox(0,0)[lb]{\smash{{\SetFigFont{12}{14.4}{\rmdefault}{\mddefault}{\updefault}{\color[rgb]{0,0,0}$\vdots$}%
}}}}
\put(6031,-2311){\makebox(0,0)[lb]{\smash{{\SetFigFont{12}{14.4}{\rmdefault}{\mddefault}{\updefault}{\color[rgb]{0,0,0}$\beta_{3,\up_1+2}$}%
}}}}
\put(5401,-2311){\makebox(0,0)[lb]{\smash{{\SetFigFont{12}{14.4}{\rmdefault}{\mddefault}{\updefault}{\color[rgb]{0,0,0}$\beta_{2,\up_1+1}$}%
}}}}
\put(4771,-2311){\makebox(0,0)[lb]{\smash{{\SetFigFont{12}{14.4}{\rmdefault}{\mddefault}{\updefault}{\color[rgb]{0,0,0}$\beta_{1,\up_1}$}%
}}}}
\put(4141,-1681){\makebox(0,0)[lb]{\smash{{\SetFigFont{12}{14.4}{\rmdefault}{\mddefault}{\updefault}{\color[rgb]{0,0,0}$\beta_{0,1}$}%
}}}}
\end{picture}%